\documentclass{amsart}
\usepackage{tikz}
\usepackage{xcolor}
\usepackage{amssymb,latexsym,amsmath,extarrows}
\usepackage{graphicx,mathrsfs,comment}
\usepackage{hyperref,url}

\newcommand{\Z}{\mathbb Z}

\renewcommand{\d}{\delta}

\renewcommand{\P}{\mathbb{P}}

\newcommand{\mc}{\mathcal}
\newcommand{\mb}{\mathbb}

\def\set4{\mathcal I}
\def\tup14{(1,2,3,4)}

\makeatletter
\newcommand\@avprod[2]{%
  {\sbox0{$\m@th#1\prod$}%
   \vphantom{\usebox0}%
   \ooalign{%
     \hidewidth
     \smash{\vrule height\dimexpr\ht0+1pt\relax depth\dimexpr\dp0+1pt\relax}%
     \hidewidth\cr
     $\m@th#1\prod$\cr
   }%
  }%
}
\newcommand{\avprod}{\mathop{\mathpalette\@avprod\relax}\displaylimits}

\usepackage{amstext}
\usepackage{bbm}

\numberwithin{equation}{section}

\newtheorem{theorem}{Theorem}[section]
\newtheorem{lemma}[theorem]{Lemma}

\newtheorem{proposition}[theorem]{Proposition}
\newtheorem{remark}[theorem]{Remark}

\newtheorem{definition}[theorem]{Definition}
\newtheorem{corollary}[theorem]{Corollary}

\newtheorem{conjecture}[theorem]{Conjecture}

\usepackage{esint}

\usepackage{enumitem}

\newcommand{\e}{\epsilon}

\newcommand{\mR}{\mathcal{R}}

\def\dm{D_{\mathrm{mul}}}

\newcommand{\R}{\mathbb{R}}

\newcommand{\supp}{\textup{supp}}

\begin{document}

\title[]{$\ell^2$ decoupling theorem for  surfaces in $\mathbb{R}^3$}

\author{Larry Guth} \address{Larry Guth\\  Deparment of Mathematics, Massachusetts Institute of Technology, USA}\email{lguth@math.mit.edu}

\author{Dominique Maldague}\address{ Dominique Maldague\\ Deparment of Mathematics, Massachusetts Institute of Technology, USA} \email{dmal@mit.edu}

\author{Changkeun Oh} \address{ Changkeun Oh\\  Deparment of Mathematics, Massachusetts Institute of Technology, 
USA, and Department of Mathematical Sciences and RIM, Seoul National University, Republic of Korea}\email{changkeun.math@gmail.com}

\maketitle

\begin{abstract}
We identify a new way to divide the $\d$-neighborhood of surfaces $\mc{M}\subset\R^3$ into a finitely-overlapping collection of rectangular boxes $S$. We obtain a sharp $(\ell^2,L^p)$ decoupling estimate using this decomposition, for the sharp range of exponents $2\le p\le 4$.  Our decoupling inequality leads to new exponential sum estimates where the frequencies lie on surfaces which do not contain a line. 


\end{abstract}




\section{Introduction}

Consider the manifold
\begin{equation}
    \mathcal{M}_{\phi}:= \{(\xi_1,\xi_2,\phi(\xi_1,\xi_2)): \xi_1,\xi_2 \in [0,1] \},
\end{equation}
where $\phi:\R^2\to\R$ is a smooth function. Denote by $N_{\delta}(\mathcal{M}_{\phi})\subset\R^3$ the $\delta$-neighborhood of the manifold $\mathcal{M}_{\phi}$.
Our main theorem is an $\ell^2$ decoupling theorem for $\mathcal{M}_{\phi}$ into $(\phi,\delta)$-flat sets. This type of theorem was introduced by 
\cite{MR4175746}.

\begin{definition}\label{1011.def11}
    Let $\phi:\R^2 \rightarrow \R$ be smooth. We say that $S \subset [0,1]^2$ is $(\phi,\delta)$-flat if
    \begin{equation}\label{24.01.11.12}
        \sup_{u,v \in S}|\phi(u)-\phi(v)- \nabla \phi(u) \cdot (u-v)| \leq \delta.
    \end{equation}
\end{definition}

For a measurable function $f:\mathbb{R}^3 \rightarrow \mathbb{C}$ and a measurable set $S \subset \mathbb{R}^2$, denote by $f_S$ a smooth Fourier restriction of $f$ on the set $S \times \mathbb{R}$. We refer to Definition \ref{1011.def12} for the rigorous definition.

\begin{theorem}\label{1005.03}
    Let $\phi:\mathbb{R}^2 \rightarrow \mathbb{R}$ be a smooth function.
    Fix $\epsilon>0$. Then there exists a sufficiently large number $A$ depending on  $\epsilon$ and $\phi$ satisfying the following. 

 For any $\delta>0$, there exists a collection $\mathcal{S}_{\delta}$ of finitely overlapping sets $S$ such that
    \begin{enumerate}
        \item the overlapping number is $O(\delta^{-\epsilon})$ in the sense that
        \begin{equation}\label{12.19.14}
            \sum_{S \in \mathcal{S}_{\delta}} \chi_S \leq C_{\epsilon,\phi}  \delta^{-\epsilon},
        \end{equation}

        \item $S$ is $(\phi,A\delta)$-flat and of the form of a parallelogram,

        \item for $2 \leq p \leq 4$ we have
        \begin{equation}\label{23.10.20}
            \|f\|_{L^p} \leq C_{\epsilon,\phi}\delta^{-\epsilon} \Big( \sum_{S \in \mathcal{S}_{\delta} }\|f_S\|_{L^p}^2 \Big)^{1/2}
        \end{equation}
        for all functions $f$ whose Fourier supports are in $N_{\delta}(\mathcal{M}_{\phi})$.
    \end{enumerate}
\end{theorem}

 One important aspect of our theorem is that the collection $\mathcal{S}_{\delta}$ is not a strict partition of $[0,1]^2$ since the sets $S\in\mc{S}_\delta$ are $O(\delta^{-\epsilon})$-overlapping
 . The decoupling inequality \eqref{23.10.20} for the hyperbolic paraboloid is not true if a family $\mathcal{S}_{\delta}$  is a partition of $[0,1]^2$. In other words, the $O( \delta^{-\epsilon})$-overlapping property is necessary for the theorem to hold true (see Appendix \ref{12.19.appendix}).
The situation is different for curves $(\xi,\phi(\xi))$ in $\R^2$, where Yang proved a version of Theorem \ref{1005.03} using a strict partition \cite{MR4310157}. We realized that to prove an $\ell^2$ decoupling inequality of the form \eqref{23.10.20}, it is crucial to understand the number of possible choices of $(\phi,\delta)$-flat sets. For a curve $(\xi,\phi(\xi))\in\mathbb{R}^2$ and a given point $p \in [0,1]$, there is essentially one choice of $(\phi,\delta)$-flat set containing the point $p$ with the maximal size. On the other hand, for the hyperbolic paraboloid in $\mathbb{R}^3$ and a given point $p \in [0,1]^2$, there are essentially $O(\log \delta^{-1})$ many choices of $(\phi,\delta)$-flat sets, which is the form of a rectangle with dyadic sidelengths, containing the point $p$ and with maximal size. This turns out to be the reason why a partition is not enough and we need \eqref{12.19.14}. For the hyperbolic paraboloid in $\mathbb{R}^n$ with $n \geq 4$, the possible choices of $(\phi,\delta)$-flat sets containing a point is $\sim \delta^{-A_n}$ for some $A_n>0$, which is why we do not have an analogous theorem without any $\delta^{-1}$-power loss (see Appendix \ref{highdim}).  \\

Next we describe our main application of Theorem \ref{1005.03}, which is to prove discrete restriction estimates for manifolds $\mathcal{M}_{\phi}$ not containing a line segment. For a function $F$ and a set $A$, define
\begin{equation}
    \|F\|_{L^p_{\#}(A)}^p:=\frac{1}{|A|}\|F\|_{L^p(A)}^p. 
\end{equation}
Denote by $\Lambda_{\delta}$  a collection of $\delta$-separated points in $[0,1]^2$. Define $e(t):=e^{2\pi i t}$.

\begin{corollary}\label{0722.thm01}
Let $\phi$ be a polynomial of degree $d$ with coefficients bounded by $1$.
    Suppose that the manifold $\mathcal{M}_{\phi}$ does not contain a line segment. Then for $2 \leq p \leq 4$, $\epsilon>0$, and any sequence $\{ a_{\xi} \}_{\xi \in \Lambda_{\delta}}$, we have
    \begin{equation}\label{0717.02}
        \Big\| \sum_{\xi \in \Lambda_{\delta}} a_{\xi}e\big(x \cdot (\xi,\phi(\xi)) \big) \Big\|_{L^p_{\#}(B_{\delta^{-d}})} \leq C_{\phi,\epsilon} \delta^{-\epsilon} \big(\sum_{\xi \in \Lambda_{\delta} }|a_{\xi}|^2 \big)^{\frac12}.
    \end{equation}
\end{corollary}
 An $(\ell^2,L^p)$ discrete restriction estimate is an inequality of the form \eqref{0717.02}.
    Corollary \ref{0722.thm01} is sharp in the following senses. First,
    the range of $p$ of \eqref{0717.02} is optimal for the paraboloid, though we expect that a wider range of $p$ is possible for generic polynomials $\phi$. Second, the inequality \eqref{0717.02} is false for any $p>2$ for some collection $\Lambda_{\delta}$ if the manifold $\mathcal{M}_{\phi}$ contains a line segment. 
    
    The $(\ell^2,L^p)$ discrete restriction estimate \eqref{0717.02} implies the sharp $(\ell^p,L^p)$ discrete restriction estimate. \cite{li2021decoupling} obtained a sharp $\ell^p$ decoupling inequality for smooth manifolds, but their decoupling does not seem to imply $(\ell^p,L^p)$ discrete restriction estimates.
    \\

Our second application is about the
Strichartz estimate for the nonelliptic Schrodinger equation on irrational tori. Let us introduce the conjecture. For $\alpha \in \mathbb{R} \setminus \mathbb{Q}$, $\alpha<0$,  define
$\tilde{\Delta}:=\frac{1}{2\pi}  (\partial_{11}+\alpha\partial_{22})$. Let $v$ be a solution of the partial differential equation 
\begin{equation}
    i \partial_t v- \tilde{\Delta}v=0, \;\;\; v(0,x)=f(x),
\end{equation}
where the initial data $f$ is defined on the torus $\mathbb{T}^2$. We let $e^{it \tilde{\Delta}}f(x):=v(x,t)$. The Strichartz estimate for the elliptic Schrodinger equation on irrational tori, which corresponds to the case $\alpha \in \mathbb{R} \setminus \mathbb{Q}$, $\alpha>0$, was studied in \cite{MR3692323,MR4436143}. As in the case of the elliptic Schrodinger equation on irrational tori in 
 \cite[Conjecture 1.2]{MR3692323}, we make the following conjecture.

\begin{conjecture} Let $\alpha$ be an irrational number. Then 
    \begin{equation}
        \|e^{it \tilde{\Delta}  }f\|_{L^p([0,T] \times \mathbb{T}^2)} \leq C_{\epsilon,\alpha} N^{\epsilon}\|f\|_{L^2(\mathbb{T}^2)}
        \begin{cases}
            T^{\frac1p} \text{ for } 2 \leq p \leq 4\\
            T^{\frac1p}+N^{1-\frac4p} \text{ for } 4 \leq p \leq 6\\
            
            T^{\frac1p}N^{1-\frac6p}+N^{1-\frac4p} \text{ for } 6 \leq p
                \end{cases}
    \end{equation}
    for functions $f$ whose Fourier support is in $[-N,N]^2$.
\end{conjecture}

 Corollary \ref{23.10.05.cor16} gives partial progress on the conjecture for the range $2 \leq p \leq 4$.

\begin{corollary}\label{23.10.05.cor16} Let $\alpha$ be an irrational number.
 Consider $2 \leq p \leq 4$. For $T \geq N$,
\begin{equation}
    \|e^{it \tilde{\Delta}  }f\|_{L^p([0,T] \times \mathbb{T}^2)} \leq C_{\epsilon,\alpha}N^{\epsilon}T^{\frac1p}\|f\|_{L^2(\mathbb{T}^2)}
\end{equation}
for functions $f$ whose Fourier support is in $[-N,N]^2$.
\end{corollary}

An analogous theorem of Corollary \ref{23.10.05.cor16} for the elliptic Schrodinger equation is obtained as a corollary of an $\ell^2$ decoupling theorem \cite{MR3374964} for the paraboloid, which is observed by \cite{MR3692323}. Corollar \ref{23.10.05.cor16} is not always true if $\alpha$ is not irrational. We refer to \cite{MR4196770, MR3321075} for related topics.

We prove the corollaries in Section \ref{corsec}.

\subsection{Comparison to decoupling inequalities in the literature}

Let $\mathcal{M}\subset\R^3$ be a manifold in $\R^3$ with nonvanishing Gaussian curvature and let $\mc{P}_\d(\mathcal{M})$ be a partition of $N_\d(\mathcal{M})$ into approximate $\d^{1/2}\times \d^{1/2}\times \d$ caps $\theta$.  In their foundational decoupling paper \cite{MR3374964}, Bourgain and Demeter proved that if $\mathcal{M}$ has everywhere positive Gaussian curvature and $2\le p\le 4$, then  
\begin{equation} \label{posdec} \|f\|_{L^p(\R^3)}\lesssim_\e \d^{-\e}(\sum_{\theta\in\mc{P}_\d(\mathcal{M})}\|f_\theta\|_{L^4(\R^3)}^2)^{1/2}. \end{equation}
for all $f$ with $\supp\widehat{f}\subset N_\d(\mathcal{M})$. It is observed in \cite{MR3736493} that \eqref{posdec} is false for the hyperbolic paraboloid. The counterexample comes from the geometric fact that the hyperbolic paraboloid contains a line. Based on this observation, one might expect \eqref{posdec} to hold true for manifolds avoiding a line segment.

\begin{conjecture}\label{24.01.13.conj16}
Let $\phi$ be a polynomial of degree $d$. Assume that $\mathcal{M}_{\phi}$ does not contain a line segment and has everywhere negative Gaussian curvature. Define
    \begin{equation}
        \mathcal{P}_{\delta}:=\{ [a,a+\delta^{\frac1d}] \times [b,b+\delta^{\frac1d}]: a,b\in \delta^{\frac1d} \mathbb{Z} \cap [0,1-\delta^{\frac1d}]   \}.
    \end{equation}
    Then for $2 \leq p \leq 4$
    \begin{equation}\label{24.01.13.112}
\|f\|_{L^p(\R^3)}\lesssim_\e \d^{-\e}(\sum_{\theta\in\mc{P}_\d}\|f_\theta\|_{L^p(\R^3)}^2)^{1/2}  \end{equation}
for
functions $f$ whose Fourier supports are in $N_{\delta}(\mathcal{M})$.
\end{conjecture}

 Conjecture \ref{24.01.13.conj16} implies Corollary \ref{0722.thm01} for the case that the determinant of the Hessian matrix of $\phi$ is negative. Unfortunately, it is not clear how to prove Conjecture \ref{24.01.13.conj16} using the current decoupling techniques. Instead of proving Conjecture \ref{24.01.13.conj16}, we introduce a modified version of the decoupling inequality, which is inspired by $(\phi,\delta)$-flat sets, and proved that this decoupling still implies Corollary \ref{0722.thm01}. Both Conjecture \ref{24.01.13.conj16} and Theorem \ref{1005.03} imply Corollary \ref{0722.thm01}. But neither of them does not seem to imply the other.
\\

Let us review decoupling for smooth manifolds to motivate the use of $(\phi,\delta)$-sets. Bourgain and Demeter \cite{MR3374964, MR3736493} proved that if $\mathcal{M}$ has everywhere nonzero Gaussian curvature and $2\le p\le 4$, then  
\begin{equation} 
\|f\|_{L^p(\R^3)}\lesssim_\e \d^{-\e}(\d^{-1})^{\frac{1}{2}-\frac{1}{p}}(\sum_{\theta\in\mc{P}_\d(\mathcal{M})}\|f_\theta\|_{L^p(\R^3)}^p)^{1/p}  \end{equation}
for all $f$ with $\supp\widehat{f}\subset N_\d(\mathcal{M})$. As mentioned before, if the Gaussian curvature is positive, then the above decoupling inequality may be refined to an $(\ell^2,L^p)$ decoupling (see \cite{MR3374964}), which has the form 
\begin{equation} \|f\|_{L^p(\R^3)}\lesssim_\e \d^{-\e}(\sum_{\theta\in\mc{P}_\d(\mathcal{M})}\|f_\theta\|_{L^p(\R^3)}^2)^{1/2}. \end{equation}
When $\mathcal{M}$ is the truncated paraboloid $\P^2=\{(\xi,|\xi|^2)\in\R^{3}:|\xi|\le 1\}$, $\mc{P}_\d(\P^2)$ is the coarsest partition of $N_\d(\P^2)$ so that each piece $\theta$ is essentially flat. Bourgain and Demeter \cite{MR3374964} also proved sharp decoupling estimates for the cylinder $Cyl=\{(\xi_1,\xi_2,\xi_3):\xi_1^2+\xi_2^2=1,\,\, |\xi_3|\lesssim 1\}$ and the cone $\mc{C}=(\xi_1,\xi_2,\sqrt{\xi_1^2+\xi_2^2}):\frac{1}{4}\le \xi_1^2+\xi_2^2\le 4\}$, where their $\d$-neighborhoods are partitioned into planks $\theta$ of dimensions about $1\times \d^{1/2}\times\d$. They observed that these are the coarsest $(\phi,\delta)$-flat sets. Based on this, Bourgain, Demeter, and Kemp \cite{MR4175746} introduced the notion of $(\phi,\delta)$-flat sets (Definition \ref{1011.def11}) for a general manifold $\mathcal{M}_{\phi}$, and asked if there is a partition of $N_{\delta}(\mathcal{M}_{\phi})$, denoted by $\mathcal{P}_{\delta}(\mathcal{M}_{\phi})$, so that each element is $(\phi,\delta)$-flat and for $2 \leq p \leq 4$ the following $\ell^p$ decoupling is true
\begin{equation}\label{12.28.129}
\|f\|_{L^p(\R^3)}\lesssim_\e \d^{-\e}|\mathcal{P}_{\delta}(\mathcal{M}_{\phi})|^{\frac{1}{2}-\frac{1}{p}}(\sum_{\theta\in\mc{P}_\d(\mathcal{M}_{\phi})}\|f_\theta\|_{L^p(\R^3)}^p)^{1/p}  \end{equation}
for all functions $f$ whose Fourier supports are in $N_{\delta}(\mathcal{M}_{\phi})$. This question is answered affirmatively by \cite{li2021decoupling} (see also the references therein).
\\

One may ask if 
there is a partition of $N_{\delta}(\mathcal{M}_{\phi})$ so that each element is $(\phi,\delta)$-flat and for $2 \leq p \leq 4$ the following $\ell^2$ decoupling is true
\begin{equation}\label{12.28.130}
\|f\|_{L^p(\R^3)}\lesssim_\e \d^{-\e}(\sum_{\theta\in\mc{P}_\d(\mathcal{M}_{\phi})}\|f_\theta\|_{L^p(\R^3)}^2)^{1/2}  \end{equation}
for all functions $f$ whose Fourier supports are in $N_{\delta}(\mathcal{M}_{\phi})$. Note that \eqref{12.28.130} implies \eqref{12.28.129} by H\"{o}lder's inequality. When $S$ is the truncated hyperbolic paraboloid, Bourgain and Demeter proved that
\begin{equation}
\|f\|_{L^p(\R^3)}\lesssim_\e \d^{-\e} \d^{-\frac12(\frac12-\frac1p)} (\sum_{\theta\in\mc{P}_\d(\mathcal{M}_{\phi})}\|f_\theta\|_{L^p(\R^3)}^2)^{1/2}.  \end{equation}
They also proved that the loss of $\delta$ on the right hand side is necessary up to the loss of $\delta^{-\epsilon}$, and the sharp example comes from the geometric fact that the hyperbolic paraboloid contains a line. At first, this might look like a counterexample to the question raised in \eqref{12.28.130}. However, since there are multiple choices of $(\phi,\delta)$-flat sets, it does not give a counterexample. Also, given that $\ell^2$ decoupling theorem for a curve in $\mathbb{R}^2$ is proved by \cite{MR4310157} (which is the two-dimensional version of \eqref{12.28.130}), one might still expect that \eqref{12.28.130} is true for some collection of $(\phi,\delta)$-flat sets. In Appendix \ref{12.19.appendix}, we prove Theorem \ref{12.28.thma1}, which says that for the hyperbolic paraboloid, there is no partition, whose elements are $(\phi,\delta)$-flat, such that \eqref{12.28.130} holds true. The proof of Theorem \ref{12.28.thma1} crucially uses the assumption that a collection of $(\phi,\delta)$-flat sets is a partition. However, in all of the known applications of decouplings we are aware of, it does not matter if a collection of $(\phi,\delta)$-flat sets is a partition or $O(\delta^{-\epsilon})$-finitely overlapping. Theorem \ref{1005.03} says that $\ell^2$ decoupling theorem is true after allowing the collection $\mathcal{S}_{\delta}$ to be $O(\log \delta^{-1})$-finitely overlapping (see \eqref{12.19.14}). Our decoupling theorem is still useful in the sense that it implies Corollary \ref{0722.thm01} and \ref{23.10.05.cor16}, and Proposition \ref{impliesTS}. Also, by H\"{o}lder's inequality, our decoupling ``essentially'' implies the $\ell^p$ decoupling theorem \eqref{12.28.129} of Yang and Li.
\\

Let us finish this subsection by mentioning higher dimensional manifolds. It is conceivable that $\ell^p$ decoupling conjecture is still true in high dimensions. Let us introduce some notations to state the conjecture. For a function $\phi:\mathbb{R}^{n-1} \rightarrow \mathbb{R}$, we say $S \subset [0,1]^{n-1}$ is $(\phi,\delta)$-set if \eqref{24.01.11.12} is true.  Consider the manifold
\begin{equation}
    \mathcal{M}_{\phi}:= \{(\xi,\phi(\xi)): \xi \in [0,1]^{n-1} \}.
\end{equation}
For a measurable function $f:\mathbb{R}^n \rightarrow \mathbb{C}$ and a measurable set $S \subset \mathbb{R}^{n-1}$, denote by $f_S$ the Fourier restriction of $f$ on the set $S \times \mathbb{R}$.

\begin{conjecture}
    Let $\phi:\mathbb{R}^{n-1} \rightarrow \mathbb{R}$ be a smooth function.
    Fix $\epsilon>0$. Then there exists a sufficiently large number $A$ depending on  $\epsilon$ satisfying the following.

 For any $\delta>0$, there exists a collection $\mathcal{S}_{\delta}$ of finitely overlapping sets $S \subset [0,1]^{n-1}$ such that
    \begin{enumerate}
        \item the overlapping number is $O(\log \delta^{-1})$ in the sense that
        \begin{equation}
            \sum_{S \in \mathcal{S}_{\delta}} \chi_S \leq C_{\epsilon} \log (\delta^{-1}).
        \end{equation}

        \item $S$ is $(\phi,A\delta)$-flat.

        \item for $2 \leq p \leq 2(n+1)/(n-1)$ we have
        \begin{equation}
            \|f\|_{L^p} \leq C_{\epsilon}\delta^{-\epsilon} (\# \mathcal{S}_{\delta})^{\frac12-\frac1p} \Big( \sum_{S \in \mathcal{S}_{\delta} }\|f_S\|_{L^p}^p \Big)^{\frac1p}
        \end{equation}
        for all functions $f$ whose Fourier supports are in $N_{\delta}(\mathcal{M}_{\phi})$.
    \end{enumerate}
\end{conjecture}
Conjecture 1.7 is known for the paraboloid and hyperbolic paraboloid, but would be new for general manifolds. The $\ell^p$ decoupling inequality for manifolds in $\R^2$ and $\R^3$ are proved by \cite{MR4310157} and \cite{li2021decoupling}. It remains open for higher dimensions. One main ingredient of $\ell^p$ decoupling for manifolds in $\R^3$ is $\ell^2$ decoupling for a curve in $\R^2$. Since we now have $\ell^2$ decoupling for manifolds in $\mathbb{R}^3$, it might be possible to prove $\ell^p$ decoupling for manifolds in $\R^4$. 


\subsection{Example: Hyperbolic paraboloid}

Our theorem is even new for the hyperbolic paraboloid. Let us state this special case.

\begin{theorem}\label{10.07.thm24}
Consider $\phi(\xi_1,\xi_2)=\xi_1\xi_2$. For $2 \leq p \leq 4$,
\begin{equation}
    \|f\|_{L^4} \leq C_{\epsilon}\delta^{-\epsilon} \big( \sum_{m=-2^{-1}\log \delta^{-1} }^{2^{-1}\log \delta^{-1}} \sum_{S \in \mathcal{R}_{ 2^m \delta^{1/2}, 2^{-m}\delta^{1/2},0 } }\|f_{S}\|_{L^4}^2 \big)^{\frac12}
\end{equation}
for all functions $f$ whose Fourier support is in $N_{\delta}(\mathcal{M}_{\phi})$. 
\end{theorem}

The set $\mathcal{R}_{ 2^m \delta^{1/2}, 2^{-m}\delta^{1/2},0 }$ is defined in Subsection \ref{0913.subsec21}.

While the $\ell^p$ decoupling for the hyperbolic paraboloid does not imply the Tomas-Stein theorem for the manifold, our theorem implies the Tomas-Stein theorem for the manifold. As an application of Theorem \ref{10.07.thm24} to exponential sum estimates, we will prove Corollary \ref{23.10.05.cor16} in Section \ref{corsec}.

Theorem \ref{10.07.thm24} itself does not imply $\ell^p$ decoupling by \cite{MR3736493}. However, in the proof of the theorem, we proved a slightly stronger inequality, which is a mixture of $\ell^2$ and $\ell^p$ norms for all the intermediate scales. We do not state this as a theorem as it is very involved. But this inequality likely implies the $\ell^p$ decoupling by \cite{MR3736493}.

We state the following version of the Stein-Tomas restriction theorem, which is recorded in Theorem 1.16 and Proposition 1.27 of \cite{MR3971577}.  Let $\mb{H}$ be the truncated hyperbolic paraboloid.
\begin{theorem}[Stein-Tomas \cite{MR0358216,MR0864375}]\label{ST} For $4\le p\le \infty$, we have  
\[ \|\widehat{F}\|_{L^p(B_{\d^{-1}})}\lesssim \d^{\frac{1}{2}}  \|F\|_{L^2(N_{\d}(\mb{H}))} \]
for each $0<\d<1$, each ball $B_{\d^{-1}}\subset\R^3$, and each $F\in L^2(N_{\d}(\mb{H}))$ supported in $N_{\d}(\mb{H})$. 
\end{theorem}

\begin{proposition}\label{impliesTS}
    Theorem \ref{10.07.thm24} implies Theorem \ref{ST} up to an $\d^{-\e}$ loss. 
\end{proposition} 

\begin{proof}[Proof of Proposition \ref{impliesTS}] Theorem \ref{ST} with $p=\infty$ is trivial, so it suffices to prove the $p=4$ case and invoke interpolation for the remaining exponents $p\in(4,\infty)$. Apply Theorem \ref{10.07.thm24} with $p=4$ to obtain the estimate
\begin{equation}\label{mainthm} \|\widehat{F}\|_{L^4(B_{\d^{-1}})}\lesssim_\e \d^{-\e}\Big(\sum_{m=-2^{-1}\log \d^{-1}}^{2^{-1}\log\d^{-1}}\sum_{S\in\mc{R}_{2^m\d^{1/2},2^{-m}\d^{1/2},0}} \|(\widehat{F})_S\|_{L^4(\R^3)}^2\Big)^{\frac{1}{2}} . \end{equation}
Note that $(\widehat{F})_S=\widehat{F\chi_S}$. Using H\"{o}lder's inequality, we have for each $S\in\mc{R}_{2^m\d^{1/2},2^{-m}\d^{1/2},0}$ that
\begin{align*}
\|\widehat{F\chi_S}\|_{L^4(\R^3)}&\le \|\widehat{F\chi_S}\|_{L^\infty(\R^3)}^{\frac{1}{2}}\|\widehat{F\chi_S}\|_{L^2(\R^3)}^{\frac{1}{2}} \\
    &\le  \|{F}\|_{L^1(N_\d(S))}^{\frac{1}{2}}\|F\|_{L^2(N_\d(S))}^{\frac{1}{2}}\\
    &\le  |N_\d(S)|^{\frac{1}{4}}  \|F\|_{L^2(N_\d(S))} . 
\end{align*}
Using the above inequality to bound the right hand side of \eqref{mainthm}, we have
\begin{align*}
    \|\widehat{F}\|_{L^4(B_{\d^{-1}})}&\lesssim_\e \d^{-\e}\Big(\sum_{m=-2^{-1}\log \d^{-1}}^{2^{-1}\log \d^{-1}} \sum_{S\in\mc{R}_{2^m\d^{1/2},2^{-m}\d^{1/2},0}}|N_\d(S)|^{\frac{1}{2}}\|F\|_{L^2(N_\d(S)}^2\Big)^{\frac{1}{2}}\\
    &\lesssim_\e \d^{-\e}\max_{\substack{|m|\lesssim \log \d^{-1}\\ S\in\mc{R}_{2^m\d^{1/2},2^{-m}\d^{1/2},0}}} |N_\d(S)|^{\frac{1}{4}} \|F\|_{L^2(N_\d(\mb{H}))}. 
\end{align*}
It remains to note that $|N_\d(S)|\le \d^2$ for each $S\in\mc{R}_{2^m\d^{1/2},2^{-m}\d^{1/2},0}$ and each $|m|\lesssim\log\d^{-1}$. 
\end{proof}

\subsection{Proof strategy/structure of the paper}

In Section \ref{sec2}, we prove  $\ell^2$ decoupling inequalities for perturbed hyperbolic paraboloids using a collection of $O(\log \delta^{-1})$-overlapping $(\phi,\delta)$-flat sets. Our result is new even for the hyperbolic paraboloid. The main new idea for the hyperbolic paraboloid is an iterative way of using broad-narrow analysis, which is combined with the induction on scale argument by \cite{MR3374964, MR3548534}. To prove theorem for perturbed hyperbolic paraboloids (Theorem \ref{0821.thm21}), we combine the aforementioned argument with a restriction theory for perturbed hyperbolic paraboloids developed by \cite{MR3999679, MR4418341, MR4631038,guo2023restriction}.
In Section \ref{1011.sec3}, we use the $\ell^2$ decoupling for the perturbed hyperbolic paraboloids as a black box, and prove the $\ell^2$ decoupling for general polynomials. This argument is essentially the same as that for \cite{li2021decoupling}.  In Section \ref{corsec}, we prove Corollary \ref{0722.thm01} and \ref{23.10.05.cor16} from the decoupling theorem. In Appendix \ref{12.19.appendix}, we prove that there is no $\ell^2$ decoupling for the hyperbolic paraboloid in $\mathbb{R}^3$ using a partition. This explains why it is necessary to introduce $\log(\delta^{-1})$ many partitions in Theorem \ref{1005.03}.
 In Appendix \ref{highdim}, we prove that there does not exist a collection of $(\phi,\delta)$-flat rectangles such that $\ell^2$ decoupling for the collection is true in high dimensions without any additional power. In appendix C, we prove the sharpness of Theorem \ref{1005.03}.

\subsection{Notations}

For a function $F$ and a set $A$, define
\begin{equation}
    \|F\|_{L^p_{\#}(A)}^p:=\frac{1}{|A|}\|F\|_{L^p(A)}^p. 
\end{equation}
For real numbers $a_i$, define 
\begin{equation}
    \Big|\avprod_{i=1}^3 a_{i}\Big|:=
    \prod_{i=1}^3 |a_{i}|^{\frac13}.
\end{equation}
For a ball $B_K$ of radius $K$ centered at $c(B_K)$, define
\begin{equation}\label{12.08.129}
    w_{B_K}(x):= \big(1+\big|\frac{x-c(B_K)}{K} \big|^2 \big)^{-100}.
\end{equation}
For a measurable set $A$ and a function $f$, define
\begin{equation}
    \fint_A |f(x)|\,dx :=\frac{1}{|A|} \int_A |f(x)|\,dx. 
\end{equation}
For a ball $B$, we define
\begin{equation}
    \|f\|_{L^p_{\#}(w_B)}^p:= |B|^{-1}\|f\|_{L^p(w_B)}^p.
\end{equation}

For a rectangular box $T$, we denote by $CT$ the dilation of $T$ by a factor of $C$ with respect to its centroid.

For two non-negative numbers $A_1$ and $A_2$, we write $A_1 \lesssim A_2$ to mean that there exists a constant $C$ such that $A_1 \le C A_2$. We write $A_1 \sim A_2$ if $A_1 \lesssim A_2$ and $A_2 \lesssim A_1$. We also write $A_1 \lesssim_{\epsilon} A_2$ if there exists $C_{\epsilon}$ depending on a parameter $\epsilon$ such that $A_1 \leq C_{\epsilon}A_2$.

\subsection*{Acknowledgements}

The authors would like to thank Yuqiu Fu for  discussion at the early stage of the project. LG is supported by a Simons Investigator award. DM is supported by the National Science Foundation under Award No. 2103249.   CO was supported by the NSF grant DMS-1800274 and POSCO Science Fellowship of POSCO TJ Park Foundation.

\section{Perturbed hyperbolic paraboloids}\label{sec2}

Consider the manifold $\mathcal{M}_{\phi}$ associated with
\begin{equation}\label{08.17.11}\phi(\xi_1,\xi_2):=
    \xi_1\xi_2+ a_{2,0}\xi_1^2+a_{0,2}\xi_2^2+ \sum_{3 \leq j+k \leq d }a_{j,k}\xi_1^{j}\xi_2^{k},
\end{equation}
where the coefficients $a_{j,k}$ satisfy
\begin{equation}\label{0928.22}
    |a_{j,k}| \leq 10^{-10d}
\end{equation}
for $d \geq 3$. 

When $d=2$, the manifold $\mathcal{M}_{\phi}$ is associated with $\phi(\xi_1,\xi_2):=\xi_1\xi_2$. Note that this is a special of case of \eqref{08.17.11} as we allow $a_{j,k}$ to be zero.

\begin{definition}\label{1011.def12}
Let $S$ be a rectangle in $\R^2$. Take a smooth function $\Xi_{S}:\R^2 \rightarrow \R$  such that
\begin{enumerate}
    \item the  support of $\widehat{\Xi_S}$ is contained in $2S$.

    \item $0 \leq \widehat{\Xi_S}(\xi_1,\xi_2) \leq 1$ for all $(\xi_1,\xi_2) \in \R^2$.
    
    \item $\widehat{\Xi_S}$ is greater than 1/10 on the set $S$.

    \item $\|\Xi_S\|_{L^1} \leq 1000$.
    
\end{enumerate}
Given a function $f:\mathbb{R}^3 \rightarrow \mathbb{C}$, define $f_S:\mathbb{R}^3 \rightarrow \mathbb{C}$ by
\begin{equation}
    f_S(x_1,x_2,x_3):= \int_{\mathbb{R}^2} f(x_1-y_1,x_2-y_2,x_3) \Xi_{S}(y_1,y_2)\,dy_1dy_2.
\end{equation}
Note that this is a convolution of $f$ and $\Xi_S$ for the first two variables.
\end{definition}

\begin{theorem}[Uniform $\ell^2$ decoupling for perturbed hyperbolic paraboloids]\label{0821.thm21} 

Fix $d \geq 2$ and $\epsilon>0$. Then there exists a sufficiently large number $A$ depending on $d$ and $\epsilon$ satisfying the following.

Let $\mathcal{M}_{\phi}$ be a manifold of the form \eqref{08.17.11} satisfying \eqref{0928.22}. Then for any $\delta>0$, there exists a family $\mathcal{S}_{\delta}$ of rectangles $S \subset \mathbb{R}^2$ such that 
    \begin{enumerate}
        \item the overlapping number is $O(\log \delta^{-1})$ in the sense that
        \begin{equation}
            \sum_{S \in \mathcal{S}_{\delta}} \chi_S \leq C_{d,\epsilon} \log (\delta^{-1}),
        \end{equation}

        \item every $S$ is $(\phi,A\delta)$-flat and of the form of a parallelogram,

        \item  for all $f$ whose Fourier support is in $N_{\delta}(\mathcal{M}_{\phi})$,
        \begin{equation}
            \|f\|_{L^4} \leq C_{d,\epsilon}\delta^{-\epsilon} \Big( \sum_{S \in \mathcal{S}_{\delta} }\|f_S\|_{L^4}^2 \Big)^{1/2}.
        \end{equation}
        The constant $C_{d,\epsilon}$ is independent of the choice of $\phi$.
    \end{enumerate}
\end{theorem}

The first step of the proof is to construct the family $\mathcal{S}_{\delta}$. A restriction theory for perturbed hyperbolic paraboloids is developed by \cite{MR3999679, MR4418341, MR4631038,guo2023restriction}. We use this theory, in particular, a language in \cite{MR4631038} to construct the family.

\subsection{Construction of the family \texorpdfstring{$\mathcal{S}_{\delta}$}{}}\label{0913.subsec21}

To construct a family $\mathcal{S}_{\delta}$, we will first define null vectors of the Hessian matrix of $\phi$.  Let us recall the notions in Section 3.1.1 of \cite{MR4631038}. Consider two vectors at $\xi \in [0,1]^2$ defined by
\begin{equation}\label{12.07.25}
    w_\xi:=(-A(\xi),1), \;\;\; v_\xi:=(1,-B(\xi)).
\end{equation}
Here $A$ and $B$ are given by
\begin{equation}
    A(\xi)=\frac{\phi_{22}(\xi)}{\phi_{12}(\xi)+\sqrt{|H_{\phi}(\xi)|} }, \;\;\; B(\xi)=\frac{\phi_{11}(\xi)}{\phi_{12}(\xi)+\sqrt{|H_{\phi}(\xi)|} },
\end{equation}
and
$H_{\phi}(\xi)$ is the determinant of the Hessian matrix of $\phi$. The functions $\phi_{ij}(\xi)$ is defined by $\partial_i \partial_j \phi(\xi)$. Note that by \eqref{0928.22} we have
\begin{equation}
    |A(\xi)| \leq 10^{-5d}, \;\;\; |B(\xi)| \leq 10^{-5d}.
\end{equation}
By (3.7) of \cite{MR4631038},
the vectors $w_\xi$ and $v_\xi$  have the following property.
\begin{equation}
    (w_\xi) H_{\phi}(\xi) (w_\xi)^{T}=0, \;\;\; (v_\xi) H_{\phi}(\xi)(v_\xi)^{T}=0.
\end{equation}
Moreover,
by \eqref{0928.22}, $w_{\xi}$ and $v_{\xi}$ are linearly independent.

We are now ready to define a family $\mathcal{S}_{\delta}$. Assume that $\delta^{-1}$ is a dyadic number. Let us first define $\mathcal{R}_{\delta}$. The family $\mathcal{S}_{\delta}$ will be a subset of it.
\begin{equation}
    \mathcal{R}_{\delta}:=\bigcup_{1 \leq \alpha \leq \delta^{-1/2}:\alpha \in 2^{\mathbb{Z}} }\mathcal{R}_{\delta\alpha,\alpha^{-1}}.
\end{equation} 
We need to define $\mathcal{R}_{\delta\alpha,\alpha^{-1}}$. Every element in the set will be a rectangle with dimension $\delta \alpha \times \alpha^{-1}$. Introduce a parameter $\beta$ which will be related to the angle of the long direction of the rectangle and the $\xi_1$-axis. Define
\begin{equation}\label{0830.29}
    \mathcal{R}_{\delta\alpha,\alpha^{-1}}:= \bigcup_{\beta \in \mathbb{Z}: 0 \leq \beta \leq \delta^{-1} \alpha^{-2}\pi } \mathcal{R}_{\delta\alpha,\alpha^{-1},\delta \alpha^2 \beta}.
\end{equation}
Here $\mathcal{R}_{\delta\alpha,\alpha^{-1},\delta \alpha^2 \beta}$ is a collection of the translated copies of the rectangle of dimension $\delta \alpha \times \alpha^{-1}$ with the angle $\delta \alpha^2 \beta$ with respect to the $\xi_1$-axis so that the elements of $\mathcal{R}_{\delta\alpha,\alpha^{-1},\delta \alpha^2 \beta}$ are disjoint and their union contains $[0,1]^2$. So we have defined 
\begin{equation}
    \mathcal{R}_{\delta}=\bigcup_{1 \leq \alpha \leq \delta^{-1/2}:\alpha \in 2^{\mathbb{Z}} }\bigcup_{\beta \in \mathbb{Z}: 0 \leq \beta \leq \delta^{-1} \alpha^{-2}\pi } \mathcal{R}_{\delta\alpha,\alpha^{-1},\delta \alpha^2 \beta}.
\end{equation}
The family $\mathcal{S}_{\delta}$ in Theorem \ref{0821.thm21} will be a subcollection of $\mathcal{R}_{\delta}$. 

\begin{definition}
Let $A$ be a sufficiently large constant, which will be the constant $A$ in Theorem \ref{0821.thm21}.
    We say two rectangles $R_1,R_2$ are comparable if
    \begin{equation}
        R_1 \subset 2AR_2, \;\;\; \mathrm{and} \;\;\; R_2 \subset 2AR_1.
    \end{equation}
\end{definition}

\begin{remark}\label{0930.rem23}
    Suppose that $R_1, R_2$ are comparable. Let $L$ be an affine transformation so that $L(R_1),L(R_2)$ are rectangles. Then $L(R_1), L(R_2)$ are also comparable.
\end{remark}

Let us now explain how to choose rectangles. Let $A$ be a large number.
Let $R \in \mathcal{R}_{\delta\alpha,\alpha^{-1},\delta \alpha^2 \beta}$. If the two following conditions are satisfied, we add $R$ to $\mathcal{S}_{\delta\alpha,\alpha^{-1},\delta \alpha^2 \beta}$.
\begin{enumerate}
    \item $R$ is $(\phi,A\delta)$-flat.

    \item For any point $z \in R$, we consider a rectangle centered at $z$, of dimension $\delta \alpha \times \alpha^{-1}$, with a long direction parallel to the vector $w_z$. This rectangle is comparable to $R$.
\end{enumerate}
If the two following conditions are satisfied, 
we also add $R$ to $\mathcal{S}_{\delta\alpha,\alpha^{-1},\delta \alpha^2 \beta}$.
\begin{enumerate}
    \item $R$ is $(\phi,A\delta)$-flat.

    \item For any point $z \in R$,  we consider a rectangle centered at $z$, of dimension $\delta \alpha \times \alpha^{-1}$, with a long direction parallel to the vector $v_z$. This rectangle is comparable to $R$.
\end{enumerate}
Recall that $v_z, w_z$ are defined in \eqref{12.07.25}.
After this process, we finally have 
\begin{equation}\label{12.07.212}
    \mathcal{S}_{\delta}:=\bigcup_{\alpha} \bigcup_{\beta} \mathcal{S}_{\delta\alpha,\alpha^{-1},\delta \alpha^2 \beta}.
\end{equation}

We define a partition of unity associated with  $\mathcal{R}_{\delta\alpha,\alpha^{-1},\delta \alpha^2 \beta}$. Namely, we take smooth functions $\{\Xi_R \}_{R \in \mR_{\delta\alpha,\alpha^{-1},\delta \alpha^2 \beta} }$ so that
\begin{enumerate}
    \item the  support of $\widehat{\Xi_R}$ is contained in $2R$.

    \item $0 \leq \widehat{\Xi_R}(\xi_1,\xi_2) \leq 1$ for all $(\xi_1,\xi_2) \in \R^2$.
    
    \item $\widehat{\Xi_R}$ is greater than 1/10 on the set $R$.

    \item for any $\xi \in [0,1]^2$ we have
    \begin{equation}
        \sum_{R \in \mR_{\delta\alpha,\alpha^{-1},\delta \alpha^2 \beta}}\widehat{\Xi_R}(\xi)=1.
    \end{equation}

    \item $\|\Xi_R\|_{L^1} \leq 1000$.
    
\end{enumerate}
For a function $f$, define $f_R$ by
\begin{equation}
    f_R(x_1,x_2,x_3):= \int_{\mathbb{R}^2} f(x_1-y_1,x_2-y_2,x_3) \Xi_{R}(y_1,y_2)\,dy_1dy_2.
\end{equation}

\subsection{Properties of partitions}

We have defined partitions $\mathcal{S}_{\delta}$ (see \eqref{12.07.212}). Let us study some properties of it.

\begin{proposition}\label{12.09.prop25}
    Let $A$ be a sufficiently large constant. Then
    \begin{equation}
        \mathcal{R}_{\delta^{1/2},\delta^{1/2},0} \subset \mathcal{S}_{\delta}.
    \end{equation}
\end{proposition}

This says that our partitions contain the squares with the canonical scale. As a remark, for the hyperbolic paraboloid, note that the following decoupling is false
\begin{equation}
            \|f\|_{L^4} \leq C_{\epsilon}\delta^{-\epsilon} \Big( \sum_{S \in \mathcal{R}_{\delta^{1/2},\delta^{1/2},0} }\|f_S\|_{L^4}^2 \Big)^{1/2}
        \end{equation} for some functions $f$ whose Fourier transforms are supported on the $\delta$-neighborhood of the hyperbolic paraboloid.

        \begin{proof}
            Take $R \in \mathcal{R}_{\delta^{1/2},\delta^{1/2},0}$. Since $R$ is a square, by the definition of $\mathcal{S}_{\delta}$, it suffices to prove that $R$ is $(\phi,A\delta)$-flat. By definition, we need to prove that
    \begin{equation}\label{0930.216}
        \sup_{u,v \in R}|\phi(u)-\phi(v)- \nabla \phi(u) \cdot (u-v)| \leq A\delta.
    \end{equation}
    Since $R$ is a square of side length $\delta^{1/2}$ and $\phi$ has the Hessian matrix whose determinant is bounded by two (see \eqref{08.17.11} and \eqref{0928.22}), \eqref{0930.216} follows from an application of Taylor's theorem.
        \end{proof}

\begin{proposition}[Finitely overlapping property]\label{12.09.prop26} For any $\xi \in [0,1]^2$
\begin{equation}
    \# \{ S \in \mathcal{S}_{\delta}: \xi \in S  \} \lesssim A\log (\delta^{-1}).
\end{equation}
\end{proposition}

This proposition gives a proof of the first property of the family $\mathcal{S}_{\delta}$ in Theorem \ref{0821.thm21} (recall that $A$ is a constant depending only on $d$ and $\epsilon$). Let us give a proof.

\begin{proof}
    Fix $\xi \in [0,1]^2$. Since the number of dyadic numbers $\alpha$ with $1 \leq \alpha \leq \delta^{-1/2}$ is $O(\log ({\delta^{-1}}))$, it suffices to show that given $\alpha$, 
    \begin{equation}
        \# \{ \beta \in \mathbb{Z} : \mathrm{there \; exist\; } S \in \mathcal{S}_{\delta\alpha,\alpha^{-1},\delta\alpha^2\beta} \mathrm{\; s.t. \;} \xi \in S   \} \lesssim A.
    \end{equation}
    Denote by $S_{1,\xi}$ (and $S_{2,\xi}$) the rectangle centered at $z$, of dimension $\delta \alpha \times \alpha^{-1}$, with a long direction parallel to the vector $w_{\xi}$ (and $v_{\xi}$).

    Suppose that $\xi \in S$ for some $S \in \mathcal{S}_{\delta\alpha,\alpha^{-1},\delta\alpha^2\beta}$. Then $S$ is comparable to either $S_{1,\xi}$ or $S_{2,\xi}$. 
    Without loss of generality, we may assume that $S$ is comparable to $S_{1,\xi}$. Then the angle between the long direction of $S$ and $w_\xi$  is $\lesssim A\delta \alpha^2$. Since $\{\delta \alpha^2 \beta\}_{\beta \in \mathbb{Z}}$ is $\delta \alpha^2$-separated, there are only $\sim A$ many $\beta$ satisfying the property.  This gives the proof.
\end{proof}

The next lemma is a technical lemma. This says that the angle between the longest direction of a rectangle $S$ and $\xi_1$-axis (or $\xi_2$-axis) is small. If $S$ is a square, it is not clear how to define ``the longest direction of $S$''. So we prove such a statement under the condition \eqref{12.09.222} to guarantee that $S$ looks like a rectangle quantitatively.
 
\begin{lemma}\label{12.09.lem27} Let $A$ be a sufficiently large number. Let $S \in \mathcal{S}_{\delta\alpha,\alpha^{-1},\delta\alpha^2\beta}$. Suppose that
    \begin{equation}\label{12.09.222}
    \frac{\mathrm{Length \; of \; the \; long \; direction \; of \; S}}{\mathrm{Length \; of \; the \; short\; direction \; of \; S}} \geq A^2.
\end{equation}
Then we have
\begin{equation}
    \delta \alpha^2 \beta \leq \frac{1}{10^{5d}}.
\end{equation}
\end{lemma}

\begin{proof}
    Suppose that $S \in \mathcal{S}_{\delta\alpha,\alpha^{-1},\delta\alpha^2\beta}$.  Let us follow the proof of Proposition \ref{12.09.prop26}. Fix $\xi \in S$. Without loss of generality, we may assume that $S$ is comparable to $S_{1,\xi}$. Then by Euclidean geometry, the angle between the long direction of $S$ and $w_\xi$  is $\lesssim A\delta \alpha^2$. The condition \eqref{12.09.222} says that $\delta \alpha^2 \leq A^{-2}$. So the angle is bounded by $\lesssim A^{-1}$. On the other hand, by \eqref{12.07.25} and the normalization condition \eqref{0928.22}, the angle between $w_{\xi}$ and $\xi_1$-axis is bounded by $10^{-9d}$. Since $A$ is sufficiently large, the angle between the long direction of $S$ and $\xi_1$-axis is bounded by $10^{-5d}$. Therefore, $\delta \alpha^2 \beta$, which indicates the angle, is bounded by $10^{-5d}$.
\end{proof}

\begin{definition}
Define $D(\delta)$ to be the smallest constant such that
\begin{equation}
    \|f\|_{L^4} \leq D(\delta) \big(\sum_{R \in \mathcal{S}_{\delta} }\|f_{R}\|_{L^4}^2 \big)^{\frac12}
\end{equation}
for all manifolds $\mathcal{M}_{\phi}$ of the form \eqref{08.17.11} satisfying \eqref{0928.22} and functions $f$ whose Fourer supports are in $N_{\delta}(\mathcal{M}_{\phi})$.    
\end{definition}

One of the key properties of the partitions is the parabolic rescaling lemma. This type of rescalings is very crucial in the work of \cite{MR3374964} and \cite{li2021decoupling}. We observe that the same rescaling lemma still holds true for our collections. The proof uses a basic property of perturbed hyperbolic paraboloids (for example, see Lemma 3.2 of \cite{guo2023restriction}).

\begin{lemma}[Parabolic rescaling]\label{0819.lem12} Let $\delta \leq \sigma \leq 1$.
    Let $R' \in \mathcal{S}_{\sigma}$. Then we have
    \begin{equation}
    \|f_{R'}\|_{L^4} \leq CD(\sigma^{-1} \delta) \big(\sum_{R \in \mathcal{S}_{\delta} }\|f_{R}\|_{L^4}^2 \big)^{\frac12}
\end{equation}
for all manifolds $\mathcal{M}_{\phi}$ of the form \eqref{08.17.11} satisfying \eqref{0928.22} and functions $f$ whose Fourer supports are in $N_{\delta}(\mathcal{M}_{\phi})$.
\end{lemma}

\begin{proof} 
    Suppose that $R'$ has dimension $\sigma\alpha \times \alpha^{-1}$ for some $1 \leq \alpha \leq \sigma^{-1}$. Then the Fourier support of $f_{R'}$ is contained in $2R' \times \mathbb{R}$. By doing the triangle inequality, we may assume that the Fourier support of $f_{R'}$ is contained in 
    $B \times \mathbb{R},$
    where $B$ is
    a box with dimension $\sigma\alpha \times \alpha^{-1}$. By translation and rotation, we may assume that the Fourier support of $f_{R'}$ is $[0,\alpha^{-1}] \times [0,\sigma\alpha]$. By abusing the notation, let us still call this box $R'$. 
    After the change of variables, our phase function changes. Let us denote the new phase function by
    \begin{equation}
        \phi_0(\xi_1,\xi_2):=\xi_1\xi_2+b_{2,0}\xi_1^2+b_{0,2}\xi_2^2+\sum_{3 \leq j+k \leq d}b_{j,k}\xi_1^j\xi_2^k.
    \end{equation}
    For convenience, we introduce $b_{1,1}=0$. Note that
    \begin{equation}\label{12.08.224}
        |b_{j,k}| \lesssim 1, \;\;\; 2 \leq j+k \leq d.
    \end{equation}
    Since $R'$
     is $(\phi_0,A\sigma)$-flat, by definition,
    \begin{equation}
        \sup_{u,v \in R'}|\phi_0(u)-\phi_0(v)-\nabla \phi_0(u) \cdot (u-v) | \leq A\sigma.
    \end{equation}
Take $u=(0,0)$ and $v=(t,0)$. Since $\nabla \phi_0(u)=0$, this inequality becomes
\begin{equation}
    \sup_{0 \leq t \leq \alpha^{-1}}|\sum_{j=2}^db_{j,0}t^j | \leq A\sigma.
\end{equation}
By \cite[Proposition 7.1]{li2021decoupling}, this inequality gives
\begin{equation}\label{12.08.227}
    |b_{j,0}| \lesssim \sigma \alpha^j, \;\; 2 \leq j \leq d.
\end{equation}
Similarly, we take $u=(0,0)$ and $v=(0,t)$ where $0 \leq t \leq \sigma \alpha$. Then by the same reasoning, we obtain
\begin{equation}\label{12.08.228}
    |b_{0,j}| \lesssim \sigma (\sigma\alpha)^{-j}, \;\; 2 \leq j \leq d.
\end{equation}    
    We next do rescaling. Define $L(\xi,\eta):=(\alpha^{-1} \xi, \sigma \alpha \eta)$ and
    \begin{equation}\label{0930.220}
         \widetilde{\phi}(\xi,\eta):= \sigma^{-1}\phi_0(L(\xi,\eta)).
    \end{equation}
The function $\tilde{\phi}(\xi,\eta)$ can be rewritten as
\begin{equation}\label{12.08.230}
    \xi_1\xi_2+(\sigma^{-1}\alpha^{-2}b_{2,0})\xi_1^2+(\sigma \alpha^2 b_{0,2})\xi_2^2+\sum_{3 \leq j+k \leq d}( \sigma^{k-1}\alpha^{k-j} b_{j,k})\xi_1^j\xi_2^k.
\end{equation}
To apply the induction hypothesis, we will prove
\begin{equation}\label{12.08.231}
    |b_{j,k}| \lesssim \sigma^{-k+1}\alpha^{-k+j}, \;\;\; 2 \leq j+k \leq d.
\end{equation}
The cases for $j=0$ or $k=0$ follow from \eqref{12.08.227} and \eqref{12.08.228}. Hence, we may assume that $j \geq 1$ and $k \geq 1$.  

Let us consider the subcase that $k \leq j$. By the inequality $1 \leq \alpha \leq \sigma^{-1}$, we have $1 \leq \alpha^{-k+j}$. So what we need to prove follows from
\begin{equation}
    |b_{j,k}| \lesssim \sigma^{-k+1}.
\end{equation}
Since $k \geq 1$ and $\sigma \leq 1$, this follows from \eqref{12.08.224}. Let us next consider the subcase that $k \geq j$. By the inequality $1 \leq \alpha \leq \sigma^{-1}$, we have $\sigma^{k-j} \leq \alpha^{-k+j}$. So what we need to prove follows from
\begin{equation}
    |b_{j,k}| \lesssim \sigma^{-j+1}.
\end{equation}
Since $j \geq 1$ and $\sigma \leq 1$, this follows from \eqref{12.08.224}. 

We have now proved \eqref{12.08.231}. So the function \eqref{12.08.230} can be rewritten as
\begin{equation}
    \xi_1\xi_2 + c_{2,0}\xi_1^2+c_{0,2}\xi_2^2 + \sum_{3 \leq j+k \leq d }c_{j,k}\xi_1^j\xi_2^k,
\end{equation}
where 
\begin{equation}
    |c_{j,k}| \lesssim 1, \;\;\; 2 \leq j+k \leq d.
\end{equation}
By applying an linear transformation, we may assume that $c_{2,0}=c_{0,2}=0$. Moreover, by doing some triangle inequality, we may assume that $|c_{j,k}| \leq 10^{-10d}$. By abusing the notation, let us still denote by $\tilde{\phi}$ the new phase function. By \eqref{0930.220}, and change of variables on physical variables, the $\delta$-neighborhood of $\mathcal{M}_{\phi}$ becomes the $\sigma^{-1}\delta$-neighborhood of $\mathcal{M}_{\tilde{\phi}}$.
    For functions $g$ whose Fourier supports are in $N_{\sigma^{-1} \delta
 }(\mathcal{M}_{\tilde{\phi}})$, we have
\begin{equation}
    \|g\|_{L^4} \leq D(\sigma^{-1} \delta) \big(\sum_{R'' \in \mathcal{S}_{\sigma^{-1}\delta} }\|g_{R''}\|_{L^4}^2 \big)^{\frac12}.
\end{equation}
    We next rescale back to the original variables. Given $R'' \in \mathcal{S}_{\sigma^{-1}\delta}$, by definition, we have 
\begin{equation}
            \sup_{u,v \in R''}|\widetilde{\phi}(u)-\widetilde{\phi}(v)-\nabla \widetilde{\phi}(u) \cdot (u-v) | \leq A\sigma^{-1}\delta.
\end{equation}
By \eqref{0930.220},  we can see that $L^{-1}(R'')$ is $(\phi,A\delta)$-flat. To conclude that $L^{-1}(R'') \in \mathcal{S}_{\delta}$, it remains to show the following property: For any point $z \in L^{-1}(R'')$,  we consider a rectangle centered at $z$, of dimension $\sigma \alpha \times \alpha^{-1}$, with a long direction parallel to the vector $v_z$ (let us denote by $R_z$). This rectangle is comparable to $L^{-1}(R'')$.

Let us prove the property. Fix $z \in L^{-1}(R'')$, and let $v_z$ be a vector associated with the phase $\phi$. By Remark \ref{0930.rem23}, it suffices to show that $L(R_z)$ is comparable to $R''$.
    After calculating the Hessian matrix, we see that $L(v_z)$ is comparable to the vector ${v_{L(z)}}$ associated with the phase function $\widetilde{\phi}$. Since $R''$ belongs to $\mathcal{S}_{\sigma^{-1}\delta}$, this gives the desired property.
\end{proof}

\subsection{The broad-narrow analysis}

In this subsection, we introduce a multilinear decoupling, and show that a linear decoupling constant is comparable to a multilinear decoupling constant up to epsilon loss (Theorem \ref{0908.thm28}). This framework is introduced by \cite{MR3374964}.

\begin{definition}\label{03.27.def211} Let $n(\xi)$ be the unit normal vector of $\mathcal{M}_{\phi}$ at the point $(\xi,\phi(\xi))$. We say three points $\xi_{(1)},\xi_{(2)},\xi_{(3)} \in [0,1]^2$ are $N^{-1}$-transverse if
\begin{equation}
    |n(\xi_{(1)}) \wedge n(\xi_{(2)}) \wedge n(\xi_{(3)})| \geq N^{-1}.
\end{equation}
Three squares $\tau_1,\tau_2,\tau_3 \subset [0,1]^2$ are called $N^{-1}$-transverse if every triple $\xi_{(i)} \in \tau_i$ is $N^{-1}$-transverse.
    
\end{definition}

\begin{definition}\label{10.2.def29}

Define $D_{\mathrm{mul}}(\delta,N^{-1})$ to be the smallest constant satisfying the following: for any $N^{-1}$-transverse squares $\tau_1,\tau_2,\tau_3 \in \mathcal{R}_{a,a,0}$ for $\delta^{1/2} \leq a \leq 1$, 
\begin{equation}
    \Big\|\avprod_{i=1}^3 f_{\tau_i} \Big\|_{L^4} \leq \dm(\delta,N^{-1}) \big(\sum_{R \in \mathcal{S}_{\delta} }\|f_{R}\|_{L^4}^2 \big)^{\frac12}
\end{equation}
for all manifolds $\mathcal{M}_{\phi}$ of the form \eqref{08.17.11} satisfying \eqref{0928.22} and functions $f$ whose Fourer supports are in $N_{\delta}(\mathcal{M}_{\phi})$.
    
\end{definition}

Recall that the definition of $f_{\tau_i}$ is given in the paragraph below \eqref{12.07.212}.
The main theorem of this subsection is as follows.

\begin{theorem}\label{0908.thm28} Given $\epsilon>0$, there exists a sufficiently large number $N$ depending on $\epsilon$ so that for all $\delta>0$ we have
    \begin{equation}
        D(\delta) \leq C_{\epsilon} 
\delta^{-\epsilon}\Big( \sup_{\delta \leq \delta' \leq 1}  D_{\mathrm{mul}}(\delta',N^{-1}) +1 \Big).
    \end{equation}
\end{theorem}

\begin{proof}

We do a broad-narrow analysis by \cite{MR2860188}. Let $K$ be a sufficiently large number, which will be determined later. Take a  sufficiently large constant $K_1$ such that $K$ is sufficiently large compared to $K_1$. Fix a ball $B_K$.  We write
\begin{equation}
    f= \sum_{\tau \in \mathcal{R}_{K^{-1/2},K^{-1/2},0 } }f_{\tau}.
\end{equation}
Recall the definition in \eqref{0830.29}. Each $\tau$ is a square of side length $K^{-1/2}$. Consider a collection of squares
\begin{equation}
    \mathcal{C}:=\{ \tau \in \mathcal{R}_{K^{-1/2},K^{-1/2},0 }: \|f_{\tau}\|_{L^4(B_K)} \geq K^{-100} \|f\|_{L^4(B_K)} \}.
\end{equation}

We will consider several cases.
By a pigeonholing argument, we can see that  $\mathcal{C}$ is nonempty. Take $\tau_1 \in \mathcal{C}$. Consider the first case that all the squares $\tau \in \mathcal{C}$ are in the $(K_1)^{-1}$-neighborhood of $\tau_1$. Denote by $\tau_1'$ the neighborhood. Then by the triangle inequality we have
\begin{equation}
    \|f\|_{L^4(B_K)} \leq \Big\| \sum_{\tau \in \mathcal{R}_{K^{-1/2},K^{-1/2},0}: \tau \cap \tau_1' \neq \emptyset } f_{\tau} \Big\|_{L^4(B_K)} + \sum_{\tau \notin \mathcal{C} }\|f_{\tau}\|_{L^4(B_K)}.
\end{equation}
By the definition of $\mathcal{C}$, this gives
\begin{equation}
    \|f\|_{L^4(B_K)} \leq 2\Big\| \sum_{\tau \in \mathcal{R}_{K^{-1/2},K^{-1/2},0}: \tau \cap \tau_1' \neq \emptyset } f_{\tau} \Big\|_{L^4(B_K)}.
\end{equation}
The set $\tau_1'$ depends on a choice of $B_K$. Hence, by summing over all possible squares, we have
\begin{equation}
    \|f\|_{L^4(B_K)} \lesssim  \big(\sum_{\tau' \in \mathcal{R}_{K_1^{-1},K_1^{-1},0} } \|f_{\tau'}\|_{L^4(B_K)}^4 \big)^{\frac14}.
\end{equation}

We next consider the case that there exists $\tau_2 \in \mathcal{C}$ outside of the $(K_1)^{-1}$-neighborhood of $\tau_1$. There are two subcases. Suppose that there are $\tau_3 \in \mathcal{C}$ such that $\tau_1,\tau_2,\tau_3$ are $10K^{-1}$-transverse. Then we have
\begin{equation}
    \|f\|_{L^4(B_K)} \lesssim K^{100}\avprod\big\|  f_{\tau_i}\big\|_{L^4(B_K)}.
\end{equation}
By an application of randomization argument and uncertainty principle, we have
\begin{equation}
    \|f\|_{L^4(B_K)} \lesssim K^{200}\big\|  \avprod \widetilde{f}_{\tau_i}\big\|_{L^4(w_{B_K})},
\end{equation}
where $\widetilde{f}_{\tau_i}$ is a modulation of ${f}_{\tau_i}$. The definition of $w_{B_K}$ is given in \eqref{12.08.129}. We refer to $(2.5)\--(2.8)$ of \cite{guth2023small} for the details of the argument.
\medskip

Lastly, suppose that there does not exist such $\tau_3$. Denote by $\xi_{(i)}$ the center of $\tau_i$ for $i=1,2$. By Definition \ref{03.27.def211}, all the elements $\tau \in \mathcal{C}$ are contained in the set
\begin{equation}\label{0830.226}
   Z_0:= \{\xi \in [0,1]^2:  |n(\xi_{(1)}) \wedge n(\xi_{(2)}) \wedge n(\xi)| < 100K^{-1}  \}.
\end{equation}
Here $n(\xi)$ is the unit normal vector of the manifold \eqref{08.17.11} at the point $(\xi,\phi(\xi))$. 
By the definition of $\mathcal{C}$, we can obtain
\begin{equation}
    \|f\|_{L^4(B_K)} \lesssim \Big\| \sum_{\tau \subset Z_0 }f_{\tau}\Big\|_{L^4(B_K)}.
\end{equation}
 We apply a uniform $\ell^2$ decoupling theorem for the set $Z_0$ of \cite[Theorem 3.1]{li2021decoupling}, and obtain
\begin{equation}\label{08.19.113}
    \begin{split}
        \|f\|_{L^4(B_K)} \lesssim_{\epsilon} K^{\epsilon} \big(\sum_{S} \|f_{S}\|_{L^4(w_{B_K})}^2)^{\frac12}
    \end{split}
\end{equation}
for any $\epsilon>0$.
Here $S \subset [0,1]^2$ is a rectangular box, and the collection $\{S\}_S$ has bounded overlap. Moreover,
\begin{equation}\label{0830.229}
    \bigcup_S S \subset \{ \xi \in[0,1]^2: |n(\xi_{(1)}) \wedge n(\xi_{(2)}) \wedge n(\xi)| \lesssim K^{-1} \}.
\end{equation}
The function $f_S$ is the Fourier restriction of $f$ to the rectangle $S \times \mathbb{R}$.\footnote{There is a minor technical issue. In the work of \cite{li2021decoupling}, they used a characteristic function for the cutoff function. The proof can be modified so that the cutoff is smooth. In our application, we use a smooth cutoff function.}
Define
\begin{equation}
    Z:=\{ \xi \in[0,1]^2: |n(\xi_{(1)}) \wedge n(\xi_{(2)}) \wedge n(\xi)|=0 \}.
\end{equation}
Then by the separation between $\tau_1$ and  $\tau_2$, and by the normalization conditions \eqref{08.17.11} and \eqref{0928.22}, one can see that
\begin{equation}
\eqref{0830.229} \subset N_{CK_1/K}(Z).
\end{equation}
Let us state this as a lemma.

\begin{lemma}\label{12.09.lem212}
    \begin{equation}
        \{ \xi \in[0,1]^2: |n(\xi_{(1)}) \wedge n(\xi_{(2)}) \wedge n(\xi)| \lesssim K^{-1} \}
        \subset N_{CK_1/K}(Z).
    \end{equation}
\end{lemma}
We postpone the proof of this lemma and continue the proof of Theorem \ref{0908.thm28}. By the lemma, each $S$ in \eqref{0830.229} is contained in a rectangle with dimension $1 \times CK_1K^{-1}$. Let us fix $S$. Suppose that 
\begin{equation}
    \frac{\mathrm{Length \; of \; the \; long \; direction \; of \; S}}{\mathrm{Length \; of \; the \; short\; direction \; of \; S}} \leq K_1.
\end{equation}
Then we decompose $S$ into squares of side length equal to the length of the short direction of $S$. Then the number of squares is bounded by $O(K_1)$. By Proposition \ref{12.09.prop25}, we have
\begin{equation}
    \|f_S\|_{L^4} \lesssim (K_1)^{\frac14} \big( \sum_{S' \in \mathcal{S}_{X}: S' \cap 2S \neq \emptyset } \|f_{S'}\|_{L^4}^2 \big)^{\frac12},
\end{equation}
where $K^{-2} \lesssim X \lesssim (K_1K^{-1})^2$. Since $K$ is sufficiently large compared to $K_1$, the term $(K_1)^{1/4}$ on the right hand side will not make any trouble.
\medskip

Let us consider the case that
\begin{equation}\label{12.09.256}
    \frac{\mathrm{Length \; of \; the \; long \; direction \; of \; S}}{\mathrm{Length \; of \; the \; short\; direction \; of \; S}} \geq K_1.
\end{equation}
Assume that $K_1$ is sufficiently large so that $K_1>A^2$.
By translation and rotation, we may assume that $S$ is contained in $S':=[0,1] \times [0,CK_1K^{-1}]$. By Lemma \ref{12.09.lem27}, after some change of variables, we may write the new phase function by
\begin{equation}
    \xi_1\xi_2+a\xi_1^2+b\xi_2^2+ \sum_{3 \leq j+k \leq d}b_{j,k}\xi_1^j\xi_2^k,
\end{equation}
where
\begin{equation}
    |a|+|b|+\sum_{3 \leq j+k \leq d}|b_{j,k}| \leq 10^{-d}.
\end{equation}
We will apply the following lemma. 

\begin{lemma}\label{0820.Lem24}
Let $K_2$ be a number such that $1 \leq K_2 \leq \delta^{-1/2}$.
    Let $S'=[0,1] \times [0,K_2^{-1}]$. Suppose that the phase function is \begin{equation}\label{03.27.264}
    \Phi(\xi_1,\xi_2):=\xi_1\xi_2+a\xi_1^2+b\xi_2^2+ \sum_{3 \leq j+k \leq d}b_{j,k}\xi_1^j\xi_2^k
\end{equation}
where
\begin{equation}
    |a|+|b|+\sum_{3 \leq j+k \leq d}|b_{j,k}| \leq 10^{-d}.
\end{equation}
Then for any $\epsilon>0$ we have
\begin{equation}\label{03.27.266}
    \|f_{S'}\|_{L^4(B_{K_2})} \leq C_{\epsilon,d,A}K_2^{\epsilon} \sup_{K_2 \leq L \leq (K_2)^2} \big( \sum_{T \in \mathcal{S}_{L^{-1} }: T \subset 2S' }\|f_T\|_{L^4(w_{B_{K_2}})}^2 \big)^{\frac12}
\end{equation}
for all functions $f$ whose Fourier support is in $N_{\delta}(\mathcal{M}_{\Phi})$. Here $f_{S'}$ is the Fourier restriction of $f$ to the rectangle $S' \times \mathbb{R}$.
\end{lemma}

Let us assume this lemma and finish the proof. By \eqref{08.19.113} and Lemma \ref{0820.Lem24}, we have
\begin{equation}
    \begin{split}
        \|f\|_{L^4(B_K)} \leq C_{\epsilon}  (K_1^{-1}K)^{\epsilon}K^{\epsilon} \sup_{(K_1^{-1}K) \leq L \leq (K)^2} \big( \sum_{T \in \mathcal{S}_{L^{-1}} }\|f_T\|_{L^4(w_{B_K})}^2 \big)^{\frac12}.
    \end{split}
\end{equation}
To summarize, by considering all possible cases, we have
\begin{equation}
\begin{split}
    \|f\|_{L^4} &\lesssim_{\epsilon}   
 \big(\sum_{\tau' \in \mathcal{R}_{K_1^{-1},K_1^{-1},0} } \|f_{\tau'}\|_{L^4}^4 \big)^{\frac14}
    \\&+K^{200}\Big\|  \avprod_{i=1}^3 \widetilde{f}_{\tau_i}\Big\|_{L^4}
    + K^{\epsilon} \sup_{K_1^{-1}K \leq L \leq (K)^2}(\sum_{T \in \mathcal{S}_L}\|f_T\|_{L^4}^2 )^{\frac12}
\end{split}
\end{equation}
for some transverse $\tau_i$.
By Definition \ref{10.2.def29} and Lemma \ref{0819.lem12}, we have
\begin{equation}\label{12.08.262}
    D(\delta) \lesssim_{\epsilon} D(K_1^2 \delta )+K^{200}D_{\mathrm{mul}}(\delta,10K^{-1})+ K^{\epsilon}
    \sup_{(K_1^{-1}K) \leq L \leq (K)^2} D(L \delta).
\end{equation}
We apply this inequality to the first and third terms on the right hand side of \eqref{12.08.262} repeatedly, and obtain the desired result.
\end{proof}

\begin{proof}[Proof of Lemma \ref{12.09.lem212}]
    For simplicity, we introduce $F(\xi):=n(\xi_{(1)}) \wedge n(\xi_{(2)}) \wedge n(\xi)$. By the definition of the normal vector $n(\xi)$, this can be rewritten as
    \begin{equation}
        F(\xi)= c \begin{vmatrix}
            \partial_1 \phi(\xi_{(1)}) & \partial_2 \phi(\xi_{(1)}) & -1 \\
            \partial_1 \phi(\xi_{(2)}) & \partial_2 \phi(\xi_{(2)})  & -1 \\
            \partial_1 \phi(\xi) & \partial_2 \phi(\xi) & -1
        \end{vmatrix}
    \end{equation}
    for some $|c| \sim 1$.
    After routine computations, this is equal to
    \begin{equation}
    \begin{split}
        -\partial_1\phi(\xi) \big(\partial_2 \phi(\xi_{(1)}) -\partial_2 \phi(\xi_{(2)}) \big)&+\partial_2\phi(\xi) \big(\partial_1 \phi(\xi_{(1)})
        -\partial_1 \phi(\xi_{(2)}) \big) \\&-\partial_1 \phi(\xi_{(1)})\partial_2 \phi(\xi_{(2)})+\partial_1 \phi(\xi_{(2)})\partial_2 \phi(\xi_{(1)}).
    \end{split}
    \end{equation}
    Recall that $\xi_{(1)}$ and $\xi_{(2)}$ are fixed.
    For simplicity, we write the function as
    \begin{equation}\label{12.09.267}
        F(\xi)=A\partial_1 \phi(\xi)+B \partial_2 \phi(\xi)+C.
    \end{equation}

    We first claim that
    \begin{equation}\label{12.09.268}
        \max\Big( |\partial_2 \phi(\xi_{(1)}) -\partial_2 \phi(\xi_{(2)})|,|\partial_1 \phi(\xi_{(1)}) -\partial_1 \phi(\xi_{(2)})|  \Big) \gtrsim \frac{1}{K_1}.
    \end{equation}
    The claim is equivalent to
    \begin{equation}
        \big| \nabla \phi(\xi_{(1)})-\nabla \phi(\xi_{(2)}) \big| \gtrsim \frac{1}{K_1}.
    \end{equation}
    By Taylor's theorem, note that
    \begin{equation}
        |\nabla^2 \phi(\xi_{(1)}) \cdot (\xi_{(1)}-\xi_{(2)})| \sim \big| \nabla \phi(\xi_{(1)})-\nabla \phi(\xi_{(2)}) \big|.
    \end{equation}
    Since the absolute value of eigenvalues of the Hessian matrix of $\phi$ is comparable to one, and by the separation $|\xi_{(1)}-\xi_{(2)}| \gtrsim \frac{1}{K_1}$, we have
    \begin{equation}
        |\nabla^2 \phi(\xi_{(1)}) \cdot (\xi_{(1)}-\xi_{(2)})| \gtrsim \frac{1}{K_1}.
    \end{equation}
    This completes the proof of the claim.
    \smallskip

    Let us continue the proof of Lemma \ref{12.09.lem212}. Recall \eqref{12.09.267}. Let us consider the case that $|A| \leq |B|$. The case $|A| \geq |B|$ can be dealt with similarly. To prove Lemma \ref{12.09.lem212}, let us fix $\xi \in [0,1]^2$ and suppose that
    \begin{equation}
        F(\xi)=X, \;\; |X| \lesssim K^{-1}.
    \end{equation}
    It suffices to show that there exists $\alpha \in \mathbb{R}$ such that
    \begin{itemize}
        \item $|\alpha| \lesssim K_1 K^{-1}$.
        \vspace{1mm}
        \item $F\big(\xi+(\alpha,0)\big)=0$.
    \end{itemize}
    By Taylor's theorem, since $F$ is a polynomial of degree $d-1$,
    \begin{equation}
        F\big(\xi+(\alpha,0)\big) = F(\xi)+ \partial_1 F(\xi) \alpha+ \sum_{2 \leq j \leq d-1} \frac{\partial_1^{j} F(\xi) }{j!}\alpha^j.
    \end{equation}
    By the expression \eqref{12.09.267},
    \begin{equation}
        \partial_1^jF(\xi)=A \partial_1^{j+1} \phi(\xi)+B\partial_1^j \partial_2 \phi(\xi).
    \end{equation}
    By the normalization \eqref{08.17.11} and \eqref{0928.22}, we have
    \begin{equation}
    \begin{split}
        F\big(\xi+(\alpha,0)\big) &= F(\xi)+(\partial_1 F(\xi)) \alpha + O\big(B\alpha^2 \big)
        \\&=X+(\partial_1 F(\xi)) \alpha + O\big(B\alpha^2 \big).
        \end{split}
    \end{equation}
    Note that $|\partial_1F(\xi)| \sim |B|$.
    By \eqref{12.09.268}, we have $|B| \gtrsim K_1^{-1}$ and we can find $|\alpha| \lesssim K_1K^{-1}$ such that
    \begin{equation}
        X+(\partial_1F(\xi))\alpha +O\big( B\alpha^2 \big)=0.
    \end{equation}
    This finishes the proof.
\end{proof}

\begin{proof}[Proof of Lemma \ref{0820.Lem24}]
    We use an induction on $K_2$. By taking $C_{\epsilon,d,A}$ sufficiently large, we may assume that Lemma \ref{0820.Lem24} is true for sufficiently large $K_2$  (compared to $A$).

We write our phase function \eqref{03.27.264} as
\begin{equation}
   A(\xi_1)+B(\xi_1,\xi_2):= \big(a\xi_1^2+\sum_{3 \leq j \leq d}b_{j,0}\xi_1^j \big) + (\xi_1\xi_2+b \xi_2^2 + \cdots ).
\end{equation}
On the set $(\xi_1,\xi_2) \in S'$,
\begin{equation}
    A(\xi_1)+B(\xi_1,\xi_2)=A(\xi_1)+O(K_2^{-1}).
\end{equation}
By  \cite[Proposition 7.1]{li2021decoupling},
\begin{equation}\label{1010.251}
    \sup_{0 \leq \xi_1 \leq 1}|A(\xi_1) | \sim_{d} \max \big(|a|,\max_j|b_{j,0}|\big).
\end{equation}

Let us consider the case that the left hand side of \eqref{1010.251} is smaller than $AK_2^{-1}$ for some fixed large number $A$. We would like to show that $S' \in \mathcal{S}_{(K_2)^{-1}}$. First, note that $S'$ is $(\phi,K_2^{-1})$-flat. Also, for any $z \in S'$, the angle between $\xi_1$-axis and $v_z$ is smaller than $O(AK_2^{-1})$. Hence, the rectangle, centered at $z$ of dimension $1 \times K_2^{-1}$ with a long direction parallel to the vector $v_z$, is comparable to $S'$. So we see that $S' \in \mathcal{S}_{(K_2)^{-1}}$. This already gives the desired result \eqref{03.27.266}.
\smallskip

We next consider the case that the left hand side of \eqref{1010.251} is larger than $AK_2^{-1}/2$. Write our function as
\begin{equation}
    f(x)=\int_{\R^3} \hat{f}(\xi)e^{2\pi i  x \cdot \xi }\,d\xi.
\end{equation}
We do a change of variables  $x_3 \mapsto M^{-1}x_3$ and $\xi_3 \mapsto M\xi_3$ so that the Fourier support of $f$ is contained in the $\delta M$-neighborhood of
\begin{equation}
    \{  M(A(\xi_1)+B(\xi_1,\xi_2)): \xi_1 \in [0,1],\, \xi_2 \in [0,K_2^{-1}]  \}
\end{equation}
for some $M^{-1} \gtrsim AK_2^{-1}$
so that $M\sup_{0 \leq \xi_1 \leq 1} |A(\xi_1)| \sim 1$. Note that
\begin{equation}
    M(A(\xi_1)+B(\xi_1,\xi_2))=MA(\xi_1)+O(A^{-1}).
\end{equation}
After the change of variables, the ball $B_{K_2}$ on the physical ball becomes a ball of radius $K_2/M$, which is still larger than $A$.
Then we use a uniform $\ell^2$ decoupling for a curve of \cite[Theorem 1.4]{MR4310157} by ignoring $\xi_2$-variable, and obtain
     \begin{equation}\label{0820.222}
    \|f_{S'}\|_{L^4} \leq C_{\epsilon/2}^{\mathrm{curve}} (A')^{\epsilon/2} \big( \sum_{T'}\|f_{T'}\|_{L^4}^2 \big)^{\frac12}.
\end{equation}
Here the sidelength of the longest direction of $T'$ is $(A')^{-1}$, which is smaller than or equal to $A^{-1/d}$. The set $T'$ has dimension $(A')^{-1} \times (K_2)^{-1}$. We next do a isotropic rescaling (with translation), and $T'$ becomes a rectangle with dimension $1 \times A'/K_2$. Since $\frac{A'}{K_2}>\frac{1}{K_2}$, by  applying the induction hypothesis on $K_2$, and rescaling back, the left hand side of \eqref{0820.222} is bounded by
\begin{equation}
C_{\epsilon/2}^{\mathrm{curve}}C_{\epsilon} (A')^{\epsilon/2}
    (\frac{K_2}{A'})^{\epsilon} \big( \sum_{T}\|f_T\|_{L^4}^2 \big)^{\frac12}.
\end{equation}
Recall that $A' \geq A^{1/d}$. Since $A$ is sufficiently large, we have $(A')^{-\epsilon/2} \leq (C_{\epsilon/2}^{\mathrm{curve}})^{-1}$, and we can close the induction.
\end{proof}

\subsection{Bourgain-Demeter type iteration}

We have shown that the linear decoupling constant is bounded by the multilinear decoupling constant (Theorem \ref{0908.thm28}). The next step is to bound the multilinear decoupling constant. Here is an ingredient.

\begin{lemma}[Ball inflation lemma]\label{0913.lem210} 

Let $\phi$ be a function \eqref{08.17.11} satisfying \eqref{0928.22}.
Let $p=4$, $q=8/3$, and $N>0$. Let $\tau_1,\tau_2,\tau_3$ be $N^{-1}$-transverse squares. For any $\rho>0$, $\epsilon>0$ and $x_0 \in \R^3$,
    \begin{equation}
        \begin{split}
            \Big(\fint_{B(x_0,\rho^{-2})} \Big( \avprod_{i=1}^3 \big(\sum_{J \in \mathcal{R}_{\rho, \rho, 0}: J \subset \tau_i } &\|f_J\|_{L^q_{\#}(w_{B_{\rho^{-1} }(x) })}^2 \big)^{\frac12} \Big)^p\,dx  \Big)^{\frac1p}
         \\&\leq C_{\epsilon,N} \rho^{-\epsilon} \avprod_{i=1}^3 \Big(\sum_{J \in \mathcal{R}_{\rho, \rho, 0}:J \subset \tau_i }\|f_J\|_{L^q_{\#}(w_{B_{\rho^{-2} }}(x_0))}^2  \Big)^{\frac{1}{2}}
         \end{split}
    \end{equation}
    for all functions $f$ whose Fourier supports are in $N_{\rho^2}(\mathcal{M}_{\phi})$.
\end{lemma}

The proof of the ball inflation lemma is standard. We refer to \cite[Theorem 6.6]{MR3548534} for the details (see also \cite[Theorem 9.2]{MR3592159}).
\medskip

We are ready to give a proof of Theorem \ref{0821.thm21}. Let $\Gamma$ be the smallest constant such that for every $\epsilon>0$, we have
\begin{equation}\label{0913.244}
    D(\delta) \leq C_{\epsilon}\delta^{-\Gamma-\epsilon}, \; \mathrm{for \; every \; dyadic \; \delta<1},
\end{equation}
where $C_{\epsilon}$ is a constant depending on $\epsilon$. Our goal is to prove $\Gamma=0$.
\\

We introduce some notations.
Take $q:=8/3$, and note that
\begin{equation}\label{0913.245}
    \frac{1}{q}=\frac{1/2}{2}+\frac{1/2}{4}.
\end{equation}
Define
\begin{equation}
\begin{split}
    &\tilde{A}_2(b):= \Big\|\avprod_{i=1}^3 \Big(\sum_{J \in \mathcal{R}_{\delta^{b/2},\delta^{b/2},0}: J \subset \tau_i }\|f_J\|_{L^2_{\#}(w_{B(x,\delta^{-b})} ) }^2 \Big)^{\frac12}\Big\|_{L^4_{x \in \R^3}},
    \\&
    \tilde{A}_q(b):= \Big\|\avprod_{i=1}^3 \Big(\sum_{J \in \mathcal{R}_{\delta^{b/2},\delta^{b/2},0}: J \subset \tau_i }\|f_J\|_{L^q_{\#}(w_{B(x,\delta^{-b/2})} ) }^2 \Big)^{\frac12}\Big\|_{L^4_{x \in \R^3}},
    \\&
    \tilde{A}_4(b):= \Big\|\avprod_{i=1}^3 \Big(\sum_{J \in \mathcal{R}_{\delta^{b/2},\delta^{b/2},0}: J \subset \tau_i }\|f_J\|_{L^4_{\#}(w_{B(x,\delta^{-b})} ) }^2 \Big)^{\frac12}\Big\|_{L^4_{x \in \R^3}},
\end{split}
\end{equation}
where $0<b<1$.
For $0<b<1$ and $*=2,q,4$, we let $a_*(b)$ the infimum over all exponents $a$ satisfying that
\begin{equation}\label{0913.247}
    \tilde{A}_{*}(b) \lesssim_{a,N} \delta^{-a}  \Big( \sum_{R \in \mathcal{S}_{\delta}} \|f_R\|_{L^4(\R^3)}^2  \Big)^{\frac12}
\end{equation}
for every $\delta>0$, every $N^{-1}$-trasverse squares $\tau_i$, and every choice of a function $f$. Lastly, define
\begin{equation}
    a_*:=\liminf_{b \rightarrow 0} \frac{\Gamma-a_*(b)}{b}
\end{equation}
for $*=2,q,4$.

\begin{lemma}\label{0913.lem211}
We have the following inequalities.
    \begin{enumerate}
    \item $(\mathrm{setup})$ $a_*< \infty$ \;\; for $*=2,q,4$
    \item $\mathrm{(rescaling)}$ $a_4 \geq \Gamma$
    \item $(L^2- \mathrm{ orthogonality})$ $a_2 \geq 2a_q$
        \item $\mathrm{(ball \; inflation)}$ $a_q \geq a_4/2+a_2/2$
    \end{enumerate}
\end{lemma}

Notice that by combining all the inequalities in the lemma, we obtain
\begin{equation}
    \Gamma \leq a_4 \leq 2a_q - a_2 \leq 0.
\end{equation}
Hence, in order to prove Theorem \ref{0821.thm21}, it suffices to prove Lemma \ref{0913.lem211}.

\begin{proof}[Proof of Lemma \ref{0913.lem211}]

    Let us show the item $(1)$. It suffices to show
    \begin{equation}\label{10.19.364}
        \Gamma \leq Cb+a_{*}(b)
    \end{equation}
    for some constant $C$.
    By the essentially constant property and Bernstein's inequality, we have
    \begin{equation}\label{0913.250}
        \Big\|\avprod_{i=1}^3 f_{\tau_i} \Big\|_{L^4} \lesssim \delta^{-Cb} \tilde{A}_{*}(b)
    \end{equation}
    for any $*=2,q,4$ and $0<b<1$. By the definition of $a_*(b)$ (see \eqref{0913.247}), we have
    \begin{equation}
        \Big\|\avprod_{i=1}^3 f_{\tau_i} \Big\|_{L^4} \lesssim \delta^{-Cb-a_{*}(b)}
         \Big( \sum_{R \in \mathcal{S}_{\delta}} \|f_R\|_{L^4(\R^3)}^2  \Big)^{\frac12}.
    \end{equation}
    Hence, we have
    \begin{equation}
        \dm(\delta,N^{-1}) \lesssim \delta^{-Cb-a_{*}(b)}.
    \end{equation}
    Combining this with Theorem \ref{0908.thm28} gives
    \begin{equation}
        D(\delta) \lesssim_{\epsilon} \delta^{-\epsilon}\delta^{-Cb-a_{*}(b)}.
    \end{equation}
    Since $\epsilon>0$ is arbitrary, by the definition of $\Gamma$, we obtain \eqref{10.19.364}.
    \medskip

    Let us next show the item $(2)$. By the definition of $a_4$, it follows from
    \begin{equation}\label{0913.255}
        a_4(b) \leq \Gamma \cdot (1-b).
    \end{equation}
    By H\"{o}lder's inequality and Minkowski's inequality, we have
    \begin{equation}
        \tilde{A}_4(b) \lesssim \avprod_{i=1}^3 \Big( \sum_{J \in \mathcal{R}_{\delta^{b/2},\delta^{b/2},0}: J \subset \tau_i } \|f_J\|_{L^4(\R^3) }^2  \Big)^{\frac12}.
    \end{equation}
By the parabolic rescaling lemma (Lemma \ref{0819.lem12}), this is further bounded by
\begin{equation}
\begin{split}
    &\lesssim D(\delta^{1-b})\avprod_{i=1}^3 \Big( \sum_{J \in \mathcal{S}_{\delta}: J \subset \tau_i } \|f_J\|_{L^4(\R^3) }^2  \Big)^{\frac12}
    \\& \lesssim D(\delta^{1-b}) \Big( \sum_{J \in \mathcal{S}_{\delta}} \|f_J\|_{L^4(\R^3) }^2  \Big)^{\frac12}.
\end{split}
\end{equation}
By \eqref{0913.244} and \eqref{0913.247}, we have
\begin{equation}
    \delta^{-a_4(b)} \lesssim \delta^{-\Gamma \cdot (1-b)}.
\end{equation}
This gives \eqref{0913.255}.
\medskip

Let us move onto the item (3). By the $L^2$-orthogonality and H\"{o}lder's inequality, we obtain
\begin{equation}
\begin{split}
\tilde{A}_2(b) &\lesssim 
    \Big\|\avprod_{i=1}^3 \Big(\sum_{J \in \mathcal{R}_{\delta^{b},\delta^{b},0}: J \subset \tau_i }\|f_J\|_{L^2_{\#}(w_{B(x,\delta^{-b})} ) }^2 \Big)^{\frac12}\Big\|_{L^4_{x \in \R^3}}
    \\& \lesssim 
    \Big\|\avprod_{i=1}^3 \Big(\sum_{J \in \mathcal{R}_{\delta^{b},\delta^{b},0}: J \subset \tau_i }\|f_J\|_{L^q_{\#}(w_{B(x,\delta^{-b})} ) }^2 \Big)^{\frac12}\Big\|_{L^4_{x \in \R^3}}.
\end{split}
\end{equation}
The last expression is $\tilde{A}_q(2b)$, so we have
\begin{equation}
    \delta^{-a_2(b)} \lesssim \delta^{-a_q(2b)},
\end{equation}
or equivalently,
\begin{equation}
    a_2(b) \leq a_q(2b).
\end{equation}
After some computations, this gives
\begin{equation}
    a_2 \geq 2a_q.
\end{equation}
This completes the proof of the item (3).
\medskip

Lastly, let us show the item (4). By the ball inflation lemma (Lemma \ref{0913.lem210}),
\begin{equation}
    \tilde{A}_q(b) \lesssim_{\epsilon} \delta^{-\epsilon} \Big\|\avprod_{i=1}^3 \Big(\sum_{J \in \mathcal{R}_{\delta^{b/2},\delta^{b/2},0}: J \subset \tau_i }\|f_J\|_{L^q_{\#}(w_{B(x,\delta^{-b})} ) }^2 \Big)^{\frac12}\Big\|_{L^4_{x \in \R^3}}.
\end{equation}
By H\"{o}lder's inequality and \eqref{0913.245},
\begin{equation}
    \tilde{A}_q(b) \lesssim_{\epsilon} \delta^{-\epsilon} \tilde{A}_2(b)^{\frac12} \tilde{A}_4(b)^{\frac12}.
\end{equation}
By the definition of $a_{*}(b)$, this gives
\begin{equation}
    \delta^{-a_q(b)} \lesssim \delta^{-\epsilon}\delta^{-\frac12(a_2(b)+a_4(b))}.
\end{equation}
After some computations, this gives
\begin{equation}
    a_q \geq a_4/2+a_2/2.
\end{equation}
This completes the proof of the item (4).
\end{proof}

\section{General manifolds}\label{1011.sec3}

We have proved the $\ell^2$ decoupling for perturbed hyperbolic paraboloids. We will use this decoupling as a black box, and prove Theorem \ref{0722.03}. This section does not contain any novelty, and we simply follow the argument of \cite{li2021decoupling} in Section 5. 
\medskip

The proof of Theorem \ref{1005.03} can be reduced to that for $\ell^2$ decoupling for polynomials (Theorem \ref{0722.03}).  We refer to Section 2.3 of \cite{li2021decoupling} for the details.

\begin{theorem}\label{0722.03}
Fix $d \geq 2$ and $\epsilon>0$. Then there exists a sufficiently large number $A$ depending on $d$ and $\epsilon$ satisfying the following.

    Let $\phi$ be a polynomial of two variables of degree $d$ with coefficients bounded by one. For any $\delta>0$, there exists a collection $\mathcal{S}_{\delta}$ of finitely overlapping sets $S$ such that
    \begin{enumerate}
        \item the overlapping number is $O( \delta^{-\epsilon})$ in the sense that
        \begin{equation}
            \sum_{S \in \mathcal{S}_{\delta}} \chi_S \leq C_{d,\epsilon} \delta^{-\epsilon}.
        \end{equation}

        \item $S$ is $(\phi,A\delta)$-flat and of the form of a parallelogram.

        \item we have
        \begin{equation}
            \|f\|_{L^4} \leq C_{d,\epsilon}\delta^{-\epsilon} \Big( \sum_{S \in \mathcal{S}_{\delta} }\|f_S\|_{L^4}^2 \Big)^{1/2}
        \end{equation}
        for all $f$ whose Fourier support is in $N_{\delta}(\mathcal{M}_{\phi})$.
        
        The constant $C_{d,\epsilon}$ is independent of the choice of $\phi$.
    \end{enumerate}
\end{theorem}

To prove Theorem \ref{0722.03}, we use Theorem \ref{0821.thm21} as a black box. Theorem \ref{0821.thm21} is a decoupling theorem  for a phase function \eqref{08.17.11} satisfying \eqref{0928.22}, but by a simple change of variables, this theorem can be generalized to that for a polynomial function $\phi$ with bounded coefficients such that
\begin{equation}
    |H_{\phi}(\xi)| \gtrsim 1
\end{equation}
for all $\xi \in [0,1]^2$. Here $H_{\phi}$ is the determinant of the Hessian matrix of $\phi$.

\subsection{Sketch of the proof of Theorem \ref{0722.03}}




Let $M$ be a sufficiently large number, which will be determined later (see the end of Section \ref{1011.sec3}). The proof is composed of three steps.

\subsection*{Step 1. Dichotomy: a curved part and a flat part}

Decompose $[0,1]^2$ according to the size of $|H_{\phi}(\xi)|$.
\begin{equation}
\begin{split}
    &S_{\mathrm{curved}}:=\{ (\xi_1,\xi_2) \in [0,1]^2 : |H_{\phi}(\xi_1,\xi_2)| > M^{-1} \},
    \\&
    S_{\mathrm{flat}}:=\{ (\xi_1,\xi_2) \in [0,1]^2 : |H_{\phi}(\xi_1,\xi_2)| < M^{-1} \}.
\end{split}
\end{equation}
Introduce a parameter $M \ll M_1$. Decompose $[0,1]^2$ into squares of side length $M_1$. Then on each square, one can see that $|H_{\phi}|>M^{-1}/2$ because $M_1^{-1} \ll M^{-1}$ and $\nabla H_{\phi}$ is bounded. By the triangle inequality,
\begin{equation}\label{0915.32}
    \|f\|_{L^4} \leq \Big\|\sum_{\substack{|\tau|=M_1^{-1}: \\ \tau \cap S_{\mathrm{curved} }\neq \emptyset } }f_{\tau}\Big\|_{L^4} +\Big\|\sum_{\substack{|\tau|=M_1^{-1}: \\ \tau \subset S_{\mathrm{flat} }  }}f_{\tau}\Big\|_{L^4}.
\end{equation}
After applying the triangle inequality,
we apply Theorem \ref{0821.thm21} (see also the discussion below Theorem \ref{0722.03}) to the first term on the right hand side, and this gives
\begin{equation}\label{0915.33}\begin{split}\Big\|\sum_{\substack{|\tau|=M_1^{-1}: \\ \tau \cap S_{\mathrm{curved} }\neq \emptyset } }f_{\tau}\Big\|_{L^4}
    &\lesssim M_1 
\big(\sum_{\substack{|\tau|=M_1^{-1}: \\ \tau \cap S_{\mathrm{curved} }\neq \emptyset } }\|f_{\tau}\|_{L^4}^2 \big)^{\frac12}
    \\&\lesssim_{\epsilon} M_1\delta^{-\epsilon}
    \Big( \sum_{S \in \mathcal{S}_{\delta} }\|f_S\|_{L^4}^2 \Big)^{1/2}.
\end{split}
\end{equation}
Since $M_1$ is a fixed number independent of $\delta$, this already gives the desired bound.
It remains to bound the second term on the right hand side of \eqref{0915.32}.

\subsection*{Step 2. Analysis of the flat part}

We next bound the term
\begin{equation}
    \Big\|\sum_{\substack{|\tau|=M_1^{-1}: \\ \tau \subset S_{\mathrm{flat} }  }}f_{\tau}\Big\|_{L^4}.
\end{equation}
By the generalized 2D uniform $\ell^2$ decoupling for polynomials (Theorem 3.1 of \cite{li2021decoupling}), we can cover $S_{\mathrm{flat}}$ by rectangles $T$ satisfying
\begin{equation}
    T \subset \{ (\xi_1,\xi_2) \in [0,1]^2: |H_{\phi}(\xi_1,\xi_2)| \lesssim M^{-1} \},
\end{equation}
and 
\begin{equation}
    \Big\|\sum_{\substack{|\tau|=M_1^{-1}: \\ \tau \subset S_{\mathrm{flat} }  }}f_{\tau}\Big\|_{L^4} \lesssim_{\epsilon} M^{\epsilon} \big( \sum_{T}\|f_T\|_{L^4}^2  \big)^{\frac12}.
\end{equation}

Let us fix $T$. Take an affine transformation $L$ such that it maps $[0,1]^2$ to $T$. Take $\phi_1:=\phi \circ L$. Then one can see that 
\begin{equation}
    |H_{\phi_1}(\xi)| \lesssim M^{-1}
\end{equation}
for all $\xi \in [0,1]^2$. By  \cite[Proposition 7.1]{li2021decoupling}, all the coefficients of $H_{\phi_1}(\xi)$ are bounded above by $\sim M^{-1}$. Hence, we can apply  \cite[Theorem 3.2]{li2021decoupling}, and obtain a rotation $\rho$ such that
\begin{equation}
    \phi_2(\xi_1,\xi_2):=\phi_1 \circ \rho (\xi_1,\xi_2)=a(\xi_1)+M^{-\alpha}b(\xi_1,\xi_2).
\end{equation}
Here $a,b$ have coefficients bounded above by $\sim 1$, and $\alpha$ is a positive number. By the uncertainty principle, the term $M^{-\alpha}B(\xi_1,\xi_2)$ is negligible on the ball of radius $M^{\alpha}$ on the physical side. Hence, we apply the 2D uniform $\ell^2$ decoupling of \cite[Theorem 4.4]{li2021decoupling} to the one variable polynomial $a(\xi_1)$, and partition $[0,1]$ into intervals $I$ so that each rectangle  is $(\phi_2,A
M^{-\alpha})$-flat.
Combining all the inequalities we have obtained so far, we have
\begin{equation}\label{0915.39}
\begin{split}
\Big\|\sum_{\substack{|\tau|=M_1^{-1}: \\ \tau \subset S_{\mathrm{flat} }  }}f_{\tau}\Big\|_{L^4} &\lesssim_{\epsilon} M^{\epsilon} \big( \sum_{T}\|f_T\|_{L^4}^2  \big)^{\frac12} 
\\&
\lesssim_{\epsilon} M^{\epsilon}   
\big( \sum_{S' \in \mathcal{S}_{M^{-\alpha}} }\|f_{S'}\|_{L^4}^2  \big)^{\frac12}. 
\end{split}
\end{equation}

\subsection*{Step 3. Iteration}

By \eqref{0915.32}, \eqref{0915.33}, and \eqref{0915.39}, we obtain
\begin{equation}\label{0915.310}
\|f\|_{L^4} \leq C_{\epsilon}
    \delta^{-\epsilon}
    \Big( \sum_{S \in \mathcal{S}_{\delta} }\|f_S\|_{L^4}^2 \Big)^{1/2}
    +C_{\epsilon}(M^{\alpha })^{\epsilon} \big( \sum_{S' \in \mathcal{S}_{M^{-\alpha}} }\|f_{S'}\|_{L^4}^2  \big)^{\frac12}.
\end{equation}
The first term is already of the desired form. To bound the second term, we fix $S'$ and do rescaling. Let $L$ be an affine transformation mapping $[0,1]^2$ to $S'$. Since $S'$ is $(\phi,AM^{-\alpha})$-flat, we have
\begin{equation}\label{12.10.314}
    \sup_{u,v \in S'}|\phi(u)-\phi(v)-\nabla \phi(u) \cdot (u-v)| \leq AM^{-\alpha}.
\end{equation}
Fix any point $u_0 \in S'$ and define
\begin{equation}
    \widetilde{\phi}(\xi):=A^{-1}M^{\alpha}\big( \phi(Lu_0)-\phi(L\xi)-\nabla \phi(Lu_0) \cdot (Lu_0 - L\xi) \big).
\end{equation}
Then by \eqref{12.10.314}, we have $\sup_{\xi \in [0,1]^2}|\widetilde{\phi}(\xi)| \leq 1$. By \cite[Proposition 7.1]{li2021decoupling}, all the coefficients of $\widetilde{\phi}$ is bounded by one. So this phase function satisfies the hypothesis of Theorem \ref{0722.03}. We apply \eqref{0915.310} to the function, and obtain
\begin{equation}\label{0915.311}
\begin{split}
    \|f_{S'}\|_{L^4}
    &\leq C_{\epsilon}
    \delta^{-\epsilon}
    \Big( \sum_{S \in \mathcal{S}_{\delta}: S \cap S' \neq \emptyset }\|f_S\|_{L^4}^2 \Big)^{1/2}
    \\&+C_{\epsilon}(M^{\alpha })^{\epsilon} \big( \sum_{S' \in \mathcal{S}_{M^{-2\alpha}}: S \cap S' \neq \emptyset }\|f_{S'}\|_{L^4}^2  \big)^{\frac12}.
\end{split}
\end{equation}
By \eqref{0915.310} and \eqref{0915.311}, we have
\begin{equation}
\begin{split}
\|f\|_{L^4} &\leq (C_{\epsilon}+C_{\epsilon}^2(M^{\alpha})^{\epsilon})
    \delta^{-\epsilon}
    \Big( \sum_{S \in \mathcal{S}_{\delta} }\|f_S\|_{L^4}^2 \Big)^{1/2}
    \\&+C_{\epsilon}^2(M^{2\alpha })^{\epsilon} \big( \sum_{S' \in \mathcal{S}_{M^{-2\alpha}} }\|f_{S'}\|_{L^4}^2  \big)^{\frac12}.
\end{split}
\end{equation}
We repeat this process $O(\log_M \delta^{-1})$-times, and  we obtain
\begin{equation}
    \|f\|_{L^4} \lesssim_{\epsilon}
    ( \delta^{-2\epsilon} + \delta^{-\log_M C_{\epsilon}} )\delta^{-\epsilon}
    \Big( \sum_{S \in \mathcal{S}_{\delta} }\|f_S\|_{L^4}^2 \Big)^{1/2}.
\end{equation}
We take $M$ sufficiently large so that $\log_M C_{\epsilon} \leq \epsilon$.
This completes the proof.

\section{Proof of Corollary \ref{0722.thm01} and \ref{23.10.05.cor16} \label{corsec}}

In this section, we prove Corollary \ref{0722.thm01} and \ref{23.10.05.cor16}. We give a remark that Corollary \ref{0722.thm01} holds true under a general condition.

\begin{remark}
    Fix $d \geq 2$.
    Let $\phi:\mathbb{R}^2 \rightarrow \mathbb{R}$ be a smooth function. Given a straight line $l$ intersecting $[0,1]^2$, we parametrize the line by $\gamma(t)$ with the unit speed. Assume that
    \begin{equation}
        \phi(\gamma(t))=a_0+a_1t+a_2t^2+\cdots+a_dt^d+E(t),
\end{equation}
where
\begin{equation}
    |E(t)| \leq c_{\phi}t^{d+1}
\end{equation}
and
\begin{equation}
    |a_2|+\cdots+|a_d| \geq C_{\phi}>0.
\end{equation}
Here, $c_{\phi}$ and $C_{\phi}$ are independent of the choice of the line $l$. Then \eqref{0717.02} is true.
\end{remark}

Corollary \ref{23.10.05.cor16} is stated using a language in a partial differential equation. By Fourier series, we can write the operator $e^{it \tilde{\Delta}}f$ as follows.
\begin{equation}
    e^{it \tilde{\Delta}}f(x)= \sum_{\xi \in [-N,N]^2 \cap \mathbb{Z}^2}\hat{f}(\xi)e\big( (x_1,x_2,x_3) \cdot (\xi_1,\xi_2,\xi_1^2-\alpha \xi_2^2) \big).
\end{equation}
Also by Parseval's identity, we have
\begin{equation}
    \|f\|_{L^2(\mathbb{T}^2)}^2 \sim \sum_{\xi \in [-N,N]^2 \cap \mathbb{Z}^2  }|\hat{f}(\xi)|^2.
\end{equation}
Hence, Corollary \ref{23.10.05.cor16} can be rephrased as follows.

\begin{corollary}\label{24.01.11.cor42}
Let $\alpha$ be irrational.
    Suppose that
    \begin{equation}
        \phi(\xi_1,\xi_2)=\xi_1^2-\xi_2^2, \;\;\; \Lambda_{\delta}=\big(\delta \mathbb{Z} \times \alpha\delta \mathbb{Z}\big) \cap [0,1]^2.
    \end{equation}
    Then for $2 \leq p \leq 4$ and $\epsilon>0$, we have
    \begin{equation}
        \Big\| \sum_{\xi \in \Lambda_{\delta}} a_{\xi}e\big(x \cdot (\xi,\phi(\xi)) \big) \Big\|_{L^p_{\#}(B_{\delta^{-3}})} \leq C_{\epsilon} \delta^{-\epsilon} \big(\sum_{\xi \in \Lambda_{\delta} }|a_{\xi}|^2 \big)^{\frac12}.
    \end{equation}
\end{corollary}

\subsection{Proof of Corollary \ref{0722.thm01}}

 Suppose that the manifold $\mathcal{M}_{\phi}$ does not contain a line. Fix $\epsilon>0$. By Theorem \ref{0722.03} with $\delta$ replaced by $\delta^d$, we have
 \begin{equation}\label{1109.51}
            \|f\|_{L^4} \leq C_{d,\epsilon}\delta^{-\epsilon} \Big( \sum_{S \in \mathcal{S}_{\delta^d} }\|f_S\|_{L^4}^2 \Big)^{1/2}
        \end{equation}
        for all $f$ whose Fourier support is in $N_{\delta^d}(\mathcal{M}_{\phi})$. Since $\phi$ does not contain a line, by a compactness argument, for any $l$, we can parametrize it as a function of $t$, and write
        \begin{equation}\label{1109.52}
            \phi(t,at+b)=a_0+a_1t+a_2t^2+\cdots+a_dt^d
        \end{equation}
        where
        \begin{equation}\label{1109.53}
            \max(|a_2|,|a_3|,\cdots,|a_d|) \sim 1.
        \end{equation}

        Suppose that $S$ is a rectangle and is $(\phi,A\delta^d)$-flat.
        We claim that the length of a long direction of $S$ is smaller or equal to $C_{\phi}\delta$ for some constant $C_{\phi}$ depending on the choice of $\phi$. We may assume that the angle between the long direction of $S$ and $\xi_1$-axis is smaller than or equal to $\pi/4$. Let $l$ be a line segment passing through the center of $S$ and parallel to the long direction of $S$ but contained in $S$. For convenience, we introduce a parametrization of the line $l$; $\gamma(t)=(t,at+b)$ where $t \in [b_0,b_1]$. To prove the claim, it suffices to show that
        \begin{equation}
            |b_1-b_0| \leq C_{\phi,A}\delta.
        \end{equation}
By the definition of $(\phi,A\delta^d)$-flat, we have
\begin{equation}
    \sup_{t \in [b_0,b_1]}|\phi(\gamma(b_0))-\phi(\gamma(t))-\nabla \phi(\gamma(b_0)) \cdot (\gamma(b_0)-\gamma(t))| \leq A\delta^d.
\end{equation}
By \eqref{1109.52} and \eqref{1109.53}, after some computations, this can be rewritten as
\begin{equation}\label{12.10.46}
    \sup_{t \in [b_0,b_1]}|\widetilde{a_1}(t-b_0)+\widetilde{a_2}(t-b_0)^2+\cdots+\widetilde{a_d}(t-b_0)^d| \leq A\delta^d,
\end{equation}
where
\begin{equation}\label{12.10.47}
    \max(|\widetilde{a_2}|,|\widetilde{a_3}|,\cdots,|\widetilde{a_d}|) \sim 1.
\end{equation}
By \cite[Proposition 7.1]{li2021decoupling},
\eqref{12.10.46} gives that 
\begin{equation}
    \max(|\widetilde{a_2}(b_0-b_1)^2|,|\widetilde{a_3}(b_0-b_1)^3|,\cdots,|\widetilde{a_d}(b_0-b_1)^d|) \lesssim A\delta^d.
\end{equation}
The condition \eqref{12.10.47} gives that
$|b_1-b_0| \lesssim \delta$, and this finishes the proof of the claim.
\smallskip

We have proved the claim. Corollary \ref{0722.thm01} will follow  by counting the number of frequencies 
 in $\Lambda_\d$ contained in a $\d^d$ neighborhood of each $S\in\mc{S}_\d$. The claim says that each $S$ contains at most $\lesssim 1$ many frequencies in $\Lambda_{\delta}$, and this gives the desired result. Let us give more details. We take $f$ to be a Fourier transform of
\begin{equation}
    \sum_{\xi \in \Lambda_{\delta}} \psi_{\xi}
\end{equation}
where $\psi_{\xi}$ is a smooth bump function supported on the ball of radius $\delta^d$ centered at the point $(\xi,\phi(\xi))$. By the claim, we have
\begin{equation}
    \Big( \sum_{S \in \mathcal{S}_{\delta^d} }\|f_S\|_{L^4}^2 \Big)^{1/2} \lesssim \Big(\sum_{\xi \in \Lambda_{\delta}} \|\widehat{\psi_{\xi}} \|_{L^4}^2 \Big)^{\frac12}.
\end{equation}
On the other hand, by the definition of $f$, we have
\begin{equation}
    \|f\|_{L^4} \sim \Big\| \sum_{\xi} \widehat{\psi_{\xi}} \Big\|_{L^4}.
\end{equation}
Note that 
\begin{equation}
 \widehat{\psi}_{\xi}(x) =\delta^{3d}e(\xi,\phi(\xi))\psi_{B(0,\delta^{-d})},  
\end{equation}
where $\psi_B$ is a smooth function essentially supported on the ball $B$. Corollary \ref{0722.thm01} follows from this, \eqref{1109.51}, \eqref{1109.52}, and \eqref{1109.53}. This completes the proof.
\\

\subsection{Proof of Corollary \ref{23.10.05.cor16}}
As we discussed, it suffices to prove Corollary \ref{24.01.11.cor42}. For simplicity, we only consider $\alpha=\sqrt{2}$. The general case can be proved identically. For $\phi(\xi_1,\xi_2)=\xi_1^2-\xi_2^2$ and $\Lambda_\delta=\delta\mathbb{Z}\times\sqrt{2}\delta\mathbb{Z}$, our goal is to show that 
\begin{equation}\label{goal}      \Big\| \sum_{\xi \in \Lambda_{\delta}} a_{\xi}e\big(x \cdot (\xi,\phi(\xi)) \big) \Big\|_{L^p_{\#}(B_{\delta^{-3}})} \leq C_{\epsilon} \delta^{-\epsilon} \big(\sum_{\xi \in \Lambda_{\delta} }|a_{\xi}|^2 \big)^{\frac12}. \end{equation}
This will follow from Theorem \ref{10.07.thm24} by counting the number of frequencies 
 in $\Lambda_\d$ contained in a $\d^3$ neighborhood of each $S\in\mc{S}_\d$. Suppose that $(m_i\d,n_i\sqrt{2}\d)\in\Lambda_\d\cap N_{\d^3}(S)$ for $i=1,2$. Then using the definition of $(\phi,\d^3)$-flat, we have
\[ |m_1^2\d^2-2n_1^2\d^2-(m_2^2\d^2-2n_2^2\d^2)-2(m_1\d,-\sqrt{2}n_1\d)\cdot((m_1-m_2)\d,\sqrt{2}(n_1-n_2)\d)|\lesssim \d^3. \]
This simplifies to 
\[ |(m_1-m_2)^2-2(n_1-n_2)^2|\lesssim \d , \]
so
\[ |[m_1-m_2-\sqrt{2}(n_1-n_2)][m_1-m_2+\sqrt{2}(n_1-n_2)]|\lesssim\d .\]
By Lemma 4 of Section 2, Chapter 2 of \cite{MR0349591}, we have $\frac{1}{b^{1+\epsilon'}}\lesssim_{\epsilon'}|a+\sqrt{2}b|$ for any $a,b\in\Z$. Therefore, if $(m_1,n_1)\not=(m_2,n_2)$, the above displayed inequality implies that $|n_1-n_2|\gtrsim_{\epsilon'} \d^{-1+\epsilon'}$. Since the elements of $\Lambda_\d\cap N_{\d^3}(S)$ are $\gtrsim \d$-separated, there are fewer than $\lesssim_{\epsilon'} \d^{-\epsilon'}$ frequencies in $\Lambda_\d\cap N_{\d^3}(S)$, as desired.

\section{Appendix: Partition is not enough}\label{12.19.appendix}

In this appendix,
we prove that there is no $\ell^2$ decoupling for the hyperbolic paraboloid using a partition. Note that in Theorem \ref{10.07.thm24} we introduce $O(\log \delta^{-1})$ many partitions to obtain $\ell^2$ decoupling for the hyperbolic paraboloid. The following theorem shows that it is necessary to introduce many partitions.

\begin{theorem}\label{12.28.thma1}
Let $\mathcal{M}_{\phi}$ be the hyperbolic paraboloid given by $\phi(\xi_1,\xi_2)=\xi_1\xi_2$. Fix $0<\epsilon<\frac{1}{100}$ and $A \geq 1$. Then the following is false: for any $\delta>0$, there exists a family $\mathcal{S}_{\delta}$ of rectangles $S \subset \mathbb{R}^2$ such that

    \begin{enumerate}
        \item the interiors of rectangles are disjoint
        \begin{equation}\label{23.12.19.11}
            \mathrm{int}(S_1) \cap \mathrm{int}(S_2) =  \emptyset, \;\;\; S_1,S_2 \in \mathcal{S}_{\delta}
        \end{equation}

        \item every $S$ is $(\phi,A\delta)$-flat

        \item  for all $f$ whose Fourier support is in $N_{\delta}(\mathcal{M}_{\phi})$,
        \begin{equation}\label{12.19.12}
            \|f\|_{L^4} \leq C_{\epsilon,A}\delta^{-\epsilon} \Big( \sum_{S \in \mathcal{S}_{\delta} }\|f_S\|_{L^4}^2 \Big)^{1/2}.
        \end{equation}
    \end{enumerate}
\end{theorem}

\begin{proof}
        Fix $\delta>0$. For simplicity, assume that $\delta^{-1}$ is a dyadic number. Suppose that such a family $\mathcal{S}_{\delta}$ exists for a contradiction. Let us first consider the partition
    \begin{equation}
    \begin{split}
        &[0,1]^2=\bigcup_{a \in \mathbb{Z} \cap [0,\delta^{-1}-1] } \Big([0,1] \times [a\delta,a\delta+\delta]\Big)=: \bigcup_a C_a,
        \\&
        [0,1]^2=\bigcup_{b \in \mathbb{Z} \cap [0,\delta^{-1}-1] } \Big([b\delta,b\delta+\delta] \times [0,1]  \Big)=: \bigcup_b D_b.
    \end{split}
    \end{equation}
    Let us fix $C_a$. Consider a collection of elements of $\mathcal{S}_{\delta}$ which have a large intersection with $C_a$
    \begin{equation}\label{12.26.14}
        \mathcal{S}_{\delta,C_a}:=\{S \in \mathcal{S}_{\delta}: |S \cap C_a| \geq \delta^{1}\delta^{\frac12-\epsilon}   \}.
    \end{equation}
    Define $\mathcal{S}_{\delta,C_b}$ similarly.  By abusing notations, let us denote by $\mathcal{S}_{\delta,a}$ and $\mathcal{S}_{\delta,b}$ the sets $\mathcal{S}_{\delta,C_a}$ and $\mathcal{S}_{\delta,D_b}$.
    Recall that elements of the set $\mathcal{S}_{\delta,a}$ are disjoint by the first condition \eqref{23.12.19.11}. 
    
    We claim that
    \begin{equation}\label{12.19.15}
        \Big|\bigcup_{S \in \mathcal{S}_{\delta,a} }S \cap C_a \,\Big| \geq \frac{99\delta}{100}, \;\;\;
        \Big|\bigcup_{S \in \mathcal{S}_{\delta,b} }S \cap D_b \,\Big| \geq \frac{99\delta}{100}.
    \end{equation}
    Let us assume this claim for a moment and finish the proof of the theorem. By taking the union over $a$ to \eqref{12.19.15}, we have
    \begin{equation}
         \Big|\bigcup_{a \in \mathbb{Z} \cap [0,\delta^{-1}-1] }\bigcup_{S \in \mathcal{S}_{\delta,a} }S \cap C_a \,\Big| \geq \frac{99}{100}.
    \end{equation}
    By pigeonholing, there exists $D_b \in \mathcal{S}_{\delta,b}$ such that
    \begin{equation}
        \Big| D_b \cap \Big( \bigcup_{a \in \mathbb{Z} \cap [0,\delta^{-1}-1] }\bigcup_{S \in \mathcal{S}_{\delta,a} }S \cap C_a \Big) \,\Big| \geq \frac{99\delta}{100}.
    \end{equation}
    This implies 
    \begin{equation}
        \Big| D_b \cap \Big( \bigcup_{a \in \mathbb{Z} \cap [0,\delta^{-1}-1] }\bigcup_{S \in \mathcal{S}_{\delta,a} }S  \Big) \,\Big| \geq \frac{99\delta}{100}.
    \end{equation}
    
    We claim that if $S \in \mathcal{S}_{\delta,a}$ for some $a$ then
    \begin{equation}\label{12.19.110}
        |D_b \cap S| \lesssim \delta^{1}\delta^{\frac12+\epsilon}.
    \end{equation}
    Note that 
    this means that $S$ does not belong to $\mathcal{S}_{\delta,D_b}$. Since the area of $D_b$ is $\delta$, the first inequality of \eqref{12.19.15} gives
\begin{equation}
        \Big|\bigcup_{S \in \mathcal{S}_{\delta,D_b} }S \cap D_b \,\Big| \leq \frac{\delta}{100}.
    \end{equation}
    This contradicts with the second inequality of \eqref{12.19.15}.
    Let us give a proof of \eqref{12.19.110}. Suppose that $S \in \mathcal{S}_{\delta,a}$. It suffices to prove that
    \begin{equation}
        S \subset (\alpha,\beta)+\Big( [0,1] \times [0,C\delta^{\frac12+\epsilon}] \Big)
    \end{equation}
    for some $\alpha,\beta \in [0,1]$. To prove this, we use the assumption that $S$ is $(\phi,A\delta)$-flat. By the definition of $\mathcal{S}_{\delta,a}$, the rectangle $S$ contains a line segment parallel to $\xi_1$-axis with length $\delta^{\frac12-\epsilon}$. Since our manifold is translation invariant, for simplicity, assume that the line is $[0,\delta^{\frac12-\epsilon}] \times \{0\}$. We will show that if $p=(p_1,p_2)$ is an element of $S$ then $|p_2| \lesssim \delta^{\frac12+\epsilon}$. This gives the proof of \eqref{12.19.110}. 
    To show the bound of $p_2$, we set $u=(0,0)$. Then by the definition of $(\phi,A\delta)$-flat (definition \ref{1011.def11}),
    \begin{equation}
        |p_1p_2| \leq A \delta.
    \end{equation}
    Similarly, we use $u=(\delta^{\frac12-\epsilon},0)$ and this gives
    \begin{equation}
        |p_1p_2-\delta^{\frac12-\epsilon}p_2| \leq A \delta.
    \end{equation}
    Combining these two gives $|p_2| \lesssim \delta^{\frac12+\epsilon}$.
\\

    It remains to prove \eqref{12.19.15}. By symmetry, let's show only the first inequality. We will use the assumption \eqref{12.19.12}. Fix $C_a$. Let us first show that there exists $S \in \mathcal{S}_{\delta}$ such that
    \begin{equation}\label{12.19.23.112}
        |S \cap C_a| \geq (\log \delta^{-1})^{-10}\delta^{1+4\epsilon}.
    \end{equation}
    Suppose that such $S$ does not exist for a contradiction. We take $f$ such that 
    \begin{enumerate}
        \item $\hat{f}$ is supported on the $\delta$-neighborhood of the set $\mathcal{M}_{\phi} \cap (C_a \times \mathbb{R})$.
        \item $\hat{f}$ is equal to one on the $\delta/2$-neighborhood of the set $\mathcal{M}_{\phi} \cap (C_a \times \mathbb{R})$.
    \end{enumerate} Here $\mathcal{M}_{\phi}$ is the hyperbolic paraboloid. Then $f$ is essentially supported on a box with dimension $\delta^{-1} \times 1 \times \delta^{-1}$ and has amplitude $\sim \delta^2$. So we have
    \begin{equation}\label{12.19.113}
        \|f\|_{L^4} \sim \delta^2 \delta^{-2/4} \sim \delta^{3/2}.
    \end{equation}
    On the other hand, by pigeonholing, there exists $j$ such that
    \begin{equation}\label{23.12.19.114}
        \Big( \sum_{S \in \mathcal{S}_{\delta} }\|f_S\|_{L^4}^2 \Big)^{1/2} \lesssim (\log \delta^{-1})^{\frac12} \Big( \sum_{S \in \mathcal{S}_{\delta}^j }\|f_S\|_{L^4}^2 \Big)^{1/2}
    \end{equation}
    where the elements $S$ of $\mathcal{S}_{\delta}^{j}$ satisfies
    \begin{equation}
        |S \cap C_a| \sim \delta 2^{-j}.
    \end{equation}
    Here, we used the fact that for large $j$ the right hand side of \eqref{23.12.19.114} is negligible (which is proved in \eqref{12.19.116}).
    Since we are assuming that there does not exist $S$ satisfying \eqref{12.19.23.112},  \eqref{23.12.19.114} is true for $j> {4\epsilon}\log_2 \delta^{-1}+10 \log_2 \log \delta^{-1}$. Note that $f_S$ is essentially constant on a box with volume $(\delta \delta 2^{-j})^{-1}$ and has amplitude $\delta \delta 2^{-j}$. Since $\mathcal{S}_{\delta}$ is a partition (see \eqref{23.12.19.11}), the cardinality of $\mathcal{S}_{\delta}^{j}$ is bounded by $2^j$. So we have
    \begin{equation}\label{12.19.116}
    \begin{split}
        \Big( \sum_{S \in \mathcal{S}_{\delta}^j }\|f_S\|_{L^4}^2 \Big)^{1/2} &\lesssim (\# \mathcal{S}_{\delta}^j)^{\frac12} \max_S \|f_S\|_{L^4}
        \\& \lesssim 2^{j/2} (\delta^2 2^{-j}) (\delta^2 2^{-j})^{-\frac14}
         \sim \delta^{\frac32}2^{-\frac{j}{4}}.
    \end{split}
    \end{equation}
    \eqref{12.19.12}, \eqref{12.19.113}, and \eqref{12.19.116} gives
    \begin{equation}
        \delta^{\frac32} \lesssim (\log \delta^{-1})^{\frac12}\delta^{\epsilon} \delta^{\frac32}2^{-\frac{j}{4}}.
    \end{equation}
    We get a contradiction from the assumption that $j>4\epsilon\log_2 \delta^{-1}+10 \log_2 \log \delta^{-1}$.  The same argument works for a general situation. Let $C \subset C_a \subset [0,1]^2$ be a convex set (polygon). Then we can show that there exists $S \in \mathcal{S}_{\delta}$ such that
    \begin{equation}\label{12.20.120}
        |S \cap C| \gtrsim(\log \delta^{-1})^{-10}\delta^{4\epsilon}|C|.
    \end{equation}
    Let us explain the proof. First, since $C$ is a convex set, we can find a rectangle $\widetilde{C} \subset C$  such that $|\widetilde{C}| \gtrsim |C|$ (for example, by Kovner–Besicovitch theorem). Since $\widetilde{C}$ is a rectangle in $C_a$, we can repeat the proof of \eqref{12.19.23.112} and find $S \in \mathcal{S}_{\delta}$ such that
    \begin{equation}
        |S \cap C| \geq |S \cap \widetilde{C}| \geq (\log \delta^{-1})^{-10}\delta^{4\epsilon}|\widetilde{C}| \gtrsim (\log \delta^{-1})^{-10}\delta^{4\epsilon}|{C}|.
    \end{equation}
    This gives \eqref{12.20.120}.
    \\

    We have shown that there exists $S_1 \in \mathcal{S}_{\delta}$ such that \eqref{12.19.23.112} holds true. In particular, $S_1 \in \mathcal{S}_{\delta,a}$ (see \eqref{12.26.14} for the definition of $\mathcal{S}_{\delta,a}$ and note that we used the assumption that $\epsilon$ is small). If we have $|S_1 \cap C_a| \geq 99\delta /100$, then this gives \eqref{12.19.15}. So suppose that
    \begin{equation}
        \delta^{1+4\epsilon} \leq |S_1 \cap C_a | \leq \frac{99}{100}\delta.
    \end{equation}
    Since $S_1$ is a rectangle, the set $C_a \setminus (S_1)^{c}$ is either a convex set or a union of two convex sets. Let us write it as
    \begin{equation}
        C_a \setminus (S_1)^c = C_{11} \cup C_{12}.
    \end{equation}
    Here $C_{11}$ and $C_{12}$ are disjoint convex sets and
    \begin{equation}
        |C_{11}|+|C_{12}| \leq |C_{a}|(1-\delta^{4\epsilon}).
    \end{equation}
    This finishes the first round of the iteration. Let us explain how to proceed.
    If $C_{11}$ satisfies $|C_{11}| \leq \delta \delta^{\frac14}$ then we leave this set. If $|C_{11}|>\delta \delta^{\frac14}$ then we apply \eqref{12.20.120} to the convex set $C_{11}$, and obtain the set $S_2 \in \mathcal{S}_{\delta}$
    such that 
    \begin{equation}
    \begin{split}
        &|S_2 \cap C_{11}|  \gtrsim (\log \delta^{-1})^{-10}\delta^{4\epsilon}|C_{11}|.
    \end{split}
\end{equation}
By the stopping time condition $|C_{11}|>\delta \delta^{\frac14}$ and the assumption that $\epsilon$ is small, we have $S_2 \in \mathcal{S}_{\delta,C_a}$ (see \eqref{12.26.14} for the definition). 
    We write $C_{11} \setminus (S_2)^c=C_{111} \cup C_{112}$. Then we have
    \begin{equation}
        |C_{111}|+|C_{112}| \leq |C_{11}|(1-\delta^{4\epsilon}),
        \end{equation}
    where $C_{11j}$ are convex sets. Repeat this process to $C_{12}$. If $|C_{12}|<\delta \delta^{\frac14}$, then we have 
    \begin{equation}
        |C_{111}|+|C_{112}|+|C_{12}| \leq |C_a|(1-\delta^{4\epsilon})^2+|C_{12}|\delta^{4\epsilon}.
    \end{equation}
    If $|C_{12}|> \delta \delta^{\frac14}$, then we have
    \begin{equation}
        |C_{111}|+|C_{112}|+|C_{113}|+|C_{114}| \leq |C_a|(1-\delta^{4\epsilon})^2.
    \end{equation}
    This finishes the second round of the iteration.
    We repeat this process $M$-times with $M \gtrsim \delta^{-6\epsilon}$. Denote by $\{C_{X}\}_{X}$ and $\{S_{Y}\}_{Y}$ a collection of convex sets and a collection of rectangles $S \in \mathcal{S}_{\delta}$ that we obtained via this process. Note that
    \begin{equation}
        \sum_{X}|C_{X}|+\sum_{Y}|S_{Y} \cap C_{a}| = |C_{a}|=\delta.
    \end{equation}
    We have
    \begin{equation}
    \begin{split}
        \sum_X|C_{X}| &= \sum_{X: |C_{X}|\geq \delta \delta^{\frac14}}|C_{X}|+\sum_{X: |C_{X}|<\delta \delta^{\frac14}}|C_{X}| 
        \\&\leq |C_a|(1-\delta^{4\epsilon})^M +\Big|\bigcup_{X:|C_X|<\delta \delta^{\frac14}}C_X \Big|\delta^{4\epsilon} \leq \frac{1}{100}|C_a|.
    \end{split}
    \end{equation}
    The last two inequalities follow from the fact that all $C_X$ are disjoint and $C_X \subset C_a$.
    This completes the proof of the claim \eqref{12.19.15}. 
 \end{proof}

\section{Appendix: Higher dimensions}\label{highdim}

In this appendix, we prove that an analogous result to Theorem \ref{1005.03} is false in higher dimensions. For convenience, let us consider only a manifold in $\mathbb{R}^4$. Higher dimensional manifolds can be proved in a similar way. Definitions \ref{1011.def11} and \ref{1011.def12} can be naturally generalized to higher dimensions. We will not state them. Recall that the critical exponent of $p$ of decoupling for the hyperbolic paraboloid in $\R^4$ is $10/3$.

\begin{theorem}\label{10.07.61}
    Consider
    $    \phi(\xi_1,\xi_2,\xi_3)=\xi_1^2+\xi_2^2-\xi_3^2.$
    Fix $0<\epsilon < \frac{1}{1000}$. Let $A$ be a constant. Then the following statement is false:

 For any $\delta>0$, there  exists a collection $\mathcal{S}_{\delta}$ of rectangular boxes $S \subset [0,1]^3$ such that
    \begin{enumerate}
        \item the overlapping number is $O(\log \delta^{-1})$ in the sense that
        \begin{equation}\label{1010.12}
            \sum_{S \in \mathcal{S}_{\delta}} \chi_S \leq C_{\epsilon} \log (\delta^{-1}).
        \end{equation}

        \item $S$ is $(\phi,A\delta)$-flat and of the form of a parallelogram.

        \item we have
        \begin{equation}\label{1010.13}
            \|f\|_{L^{10/3}} \leq C_{\epsilon}\delta^{-\epsilon} \Big( \sum_{S \in \mathcal{S}_{\delta} }\|f_S\|_{L^{10/3}}^2 \Big)^{1/2}
        \end{equation}
        for all $f$ whose Fourier support is in $N_{\delta}(\mathcal{M}_{\phi})$.
    \end{enumerate}
\end{theorem}

The same conclusion of Theorem \eqref{10.07.61} is where \eqref{1010.13} is replaced with
\begin{equation}\label{1007.63}
            \|f\|_{L^{p}} \leq C_{\epsilon}\delta^{-\epsilon} \Big( \sum_{S \in \mathcal{S}_{\delta} }\|f_S\|_{L^{p}}^2 \Big)^{1/2}
        \end{equation}
    for $2<p \leq 10/3$. \eqref{1007.66} is the only place where \eqref{1010.13} is used. One may check that \eqref{1007.63} still implies the desired result. 

\begin{proof}
    Let us prove by contradiction. Fix $0<\epsilon<10^{-10}$. Define
    \begin{equation}
        \mathcal{C}:=\{(\xi_1,\xi_2,\xi_3) \in [0,1]^3: \xi_1^2+\xi_2^2-\xi_3^2=0  \}.
    \end{equation}
    We cover it by $\sim \delta^{-1/2}$ many canonical blocks $P$ with dimension $1 \times \delta^{1/2} \times \delta$. Let us denote by $\mathcal{P}$ the collection of blocks. Note that each block is $(\phi,A\delta)$-flat for some large constant $A$. For given $P \in \mathcal{P}$, we take translated copies and cover $[0,1]^3$ so that they are disjoint. Denote by $\mathcal{P}_P$ the collection of translated copies of $P$. Note that the cardinality of $\mathcal{P}_P$ is comparable to $\delta^{-3/2}$.

    Fix $P' \in \mathcal{P}_P$ for some $P \in \mathcal{P}$. We next decompose $P'$ into parallel tubes with radius $\delta^X$, where $X$ is sufficiently large number. Let us denote by $\mathcal{P}_{P,P'}$ the collections of tubes $P''$. We claim that for given $P'' \in \mathcal{P}_{P,P'}$ there is an element $S$ of $\mathcal{S}_{\delta}$ such that the diameter of the set $S \cap P''$ is larger than or equal to $(\log \delta^{-1})\delta^{5\epsilon}$. To prove the claim, let us use the hypothesis \eqref{1010.13}. Take a function $f$ so that
    \begin{equation}
        \hat{f}(\xi):= \sum_{j=1}^{(\log \delta^{-1})\delta^{-5\epsilon}} f_j(\xi):=\sum_{j=1}^{(\log \delta^{-1})\delta^{-5\epsilon}} 1_{N_{\delta^X}( B_j)}(\xi)
    \end{equation}
    where $\{B_j\}_{j=1}^{(\log \delta^{-1})\delta^{-5\epsilon}}$ is an arithmetic progression with difference $(\log \delta^{-1})\delta^{5\epsilon}$, and $\bigcup_j B_j \subset P''$. By the direct computations, one can see that
    \begin{equation}\label{1007.66}
        \|f\|_{L^{10/3}} \gtrsim (\log \delta^{-1})^{\frac15}\delta^{-5\epsilon(\frac12-\frac{1}{10/3})} (\sum_j \|f_j\|_{L^{10/3}}^2)^{1/2}.
    \end{equation} 
    Hence, \eqref{1010.13} says that there must be an element $S$ of $\mathcal{S}_{\delta}$ containing at least two $B_j$. This means that the diameter of $S \cap P''$ is larger than or equal to $(\log \delta^{-1})\delta^{5\epsilon}$. This completes the proof of the claim.
    
    For given $P' \in \mathcal{P}_P$, denote by $\mathcal{S}_{P'}$ the elements of $S \in \mathcal{S}_{\delta}$ satisfying the property: for some $P'' \in \mathcal{P}_{P,P'}$ the diameter of $S \cap P''$ is larger than or equal to $(\log \delta^{-1})\delta^{5\epsilon}$. The claim says that $\mathcal{S}_{P'}$ is non-empty.

We next claim that if $S$ is $(\phi,A\delta)$-flat and $S \in \mathcal{S}_{\delta}$, then $S$ is contained in a  box with dimension $\delta^{\frac12-\frac{5\epsilon}{2}} \times \delta^{1-5\epsilon} \times \delta^{5\epsilon}$. Let us give  a proof. First of all, by an affine transformation, we may assume that $S$ contains a line $\{0\} \times \{0\} \times [0,\delta^{5\epsilon}]$ and our new phase function is $\widetilde{\phi}(\xi)=\xi_1^2+\xi_2\xi_3$. Let $v=(v_1,v_2,v_3) \in S$. We will show that 
\begin{equation}
    |v_1| \lesssim \delta^{\frac12-\frac{5\epsilon}{2}}, \;\;\; |v_2| \lesssim \delta^{1-5\epsilon}.
\end{equation}
By the definition of $(\phi,A\delta)$-flat set with $u=(0,0,0)$ we obtain
\begin{equation}\label{12.27.27}
    |v_1^2+v_2v_3| \leq A\delta.
\end{equation}
We next use $u=(0,0,\delta^{5\epsilon})$ and obtain
\begin{equation}
    |v_1^2+v_2v_3- \delta^{5\epsilon}v_2| \leq A\delta.
\end{equation}
These two inequalities give $|v_2| \lesssim \delta^{1-5\epsilon}$. This bound and \eqref{12.27.27} give $|v_1| \lesssim \delta^{\frac12-\frac{5\epsilon}{2}}$. This proves the claim.

 By the claim, if $P_1' \in \mathcal{P}_{P_1}$ and $P_2' \in \mathcal{P}_{P_2}$, and the angle of the longest directions of $P_1$ and $P_2$ is greater than $\delta^{\frac12-100\epsilon}$, then any two sets $S_1 \in \mathcal{S}_{P_1'}$ and $S_2 \in \mathcal{S}_{P_2'}$ are distinct.  Note also that
\begin{equation}
    \Big|\bigcup_{S' \in \mathcal{S}_{P'} } S' \cap P' \Big| \gtrsim \delta^{5\epsilon}|P'|.
\end{equation}

Let us finish the proof of the theorem. 
To get a contradiction, we need to prove that there exists $\xi \in [0,1]^3$ such that the number of $S\in \mathcal{S}_{\delta}$ containing $\xi$ is greater than the right hand side of \eqref{1010.12}. To prove this, it suffices to show that there exists $\xi\in [0,1]^3$ satisfying
\begin{equation}
    \sum_{P \in \mathcal{P}}
    \sum_{P' \in \mathcal{P}_P }
    \sum_{S' \in \mathcal{S}_{P'} }
    1_{S' \cap P' }(\xi) \gtrsim 
 \delta^{-\alpha-100\epsilon}
\end{equation}
for some number $\alpha>0$.
This follows from a pigeonholing argument with
\begin{equation}
\begin{split}
    \int_{[0,1]^3}\sum_{P \in \mathcal{P}}
    \sum_{P' \in \mathcal{P}_P }
    \sum_{S' \in \mathcal{S}_{P'} }
    1_{S' \cap P'}(\xi)  \,d\xi &=
    \sum_{P \in \mathcal{P}}
    \sum_{P' \in \mathcal{P}_P }
    \int_{[0,1]^3}
    \sum_{S' \in \mathcal{S}_{P'} }
    1_{S' \cap P'}(\xi)  \,d\xi
     \\&\gtrsim \delta^{-\frac32}\delta^{-\frac12} \delta^{5\epsilon}\delta^{\frac32} \sim \delta^{-\frac12+5\epsilon}.
\end{split}
\end{equation}
This completes the proof.
\end{proof}

\section{Appendix: Sharpness of Theorem \ref{1005.03}}

Let us discuss the sharpness of Theorem \ref{1005.03}. More precisely, we would like to discuss if the  range of $p$ for which \eqref{23.10.20} holds true is sharp. An example of a manifold that the aforementioned range of $p$ is not sharp is as follows.
\begin{equation}\label{23.10.20.19}
    \phi(\xi_1,\xi_2)=\xi_1^2.
\end{equation}
The $\ell^2$ decoupling for \eqref{23.10.20.19} is proved for $2 \leq p \leq 6$ by \cite{MR3374964}, and this range of $p$ is sharp. Motivated by this example, let us introduce the following definition.

\begin{definition}
    We say  $\phi(\xi_1,\xi_2)$ depends only on a variable if there exists an affine transformation $L$ such that
    \begin{equation}\label{23.10.20.110}
        \phi(L(\xi_1,\xi_2))=\psi(\xi_1)+a\xi_2
    \end{equation}
    for some function $\psi$ and some $a \in \mathbb{R}$. 
    Second, we say $\phi(\xi_1,\xi_2)$ does not depend on any variable if
    \begin{equation}\label{23.10.20.111}
        \phi(\xi_1,\xi_2)=a+b\xi_1+c\xi_2
    \end{equation}
    for some $a,b,c \in \mathbb{R}$. 
    Lastly, we say $\phi(\xi_1,\xi_2)$ depends on two variables if it is not a form of \eqref{23.10.20.111} and there does not exist $L$ such that \eqref{23.10.20.110} holds true.
\end{definition}

Here is a complete characterization of the $\ell^2$ decoupling theorem.

\begin{proposition}\label{10.23.18}
 Let $p \geq 2$.
Consider a polynomial  $\phi:=\phi(\xi_1,\xi_2)$. Suppose that $\mathcal{S}_{\delta}$ is a family constructed in Subsection \ref{0913.subsec21} and Section \ref{1011.sec3}.
    \begin{enumerate}
        \item Let $\phi$ depend on two variables. Then \eqref{23.10.20} is true for $2 \leq p \leq 4$. 
        
        \item Let $\phi$ depend only on a variable. Then \eqref{23.10.20} is true for $2 \leq p \leq 6$. 

        \item  Let $\phi$  not depend on any variable. Then \eqref{23.10.20} is true  for $2 \leq p <\infty$.
        
    \end{enumerate}
    The ranges of $p$ stated in items (1), (2), and (3) are sharp.
\end{proposition}

\begin{proof}

Let us first prove Item (1). Theorem \ref{1005.03} gives \eqref{23.10.20} for $2 \leq p \leq 4$, so it suffices to prove that the range is sharp. Let $H_{\phi}(\xi)$ be the Hessian matrix of $\phi$. By \cite[Theorem 3.1]{MR2095579}, $H_{\phi}(\xi)$ is not identically zero. Consider 
\begin{equation}
    Z_{\phi}:=\{\xi \in \mathbb{R}^2: \det{(H_{\phi}(\xi))}=0 \}.
\end{equation}
Since $\det H_{\phi}(\xi)$ is not identically zero, it is a zero set of a polynomial.  So there exists a square $\tau$ not intersecting $Z_{\phi}$. Since $H_{\phi}(\xi)$ is independent of the parameter $\delta$, we may assume that the size of $\tau$ is independent of $\delta$. Let us fix such $\tau$. Assume that $\delta$ is  smaller than the sidelength of $\tau$.  Let $\hat{f}$ be a smooth bump  function of the $\delta$-neighborhood of
\begin{equation}
    \{ (\xi_1,\xi_2,\phi(\xi_1,\xi_2) ) :\xi \in \tau \}.
\end{equation} 
One can see that if $\xi \in \tau$, then
\begin{equation}
    |\det{H_{\phi}(\xi)}| \gtrsim 1.
\end{equation}
So the multiplication of eigenvalues of $H_{\phi}(\xi)$ does not change the sign over $\xi \in \tau$.

Suppose that the eigenvalues of $ H_{\phi}(\xi)$ have the same sign for all $\xi \in \tau$. 
Then each $S \in \mathcal{S}_{\delta}$ intersecting $\tau$ is a square with length $\delta^{1/2} \times \delta^{1/2}$. So we have
\begin{equation}
    f=\sum_S f_S, \;\;\; |f_S(x)|=\delta^21_{T_S}
\end{equation}
where $T_S$ is a tube with dimension $\delta^{-1/2} \times \delta^{-1/2} \times \delta^{-1}$. Then
\begin{equation}
    \delta^2 \delta^{-1} \lesssim \|f\|_{L^p(B_1)} \lesssim \|f\|_{L^p(\R^3)}.
\end{equation}
On the other hand,
\begin{equation}
    \Big( \sum_S \|f_S\|_{L^p}^2 \Big)^{\frac12} \lesssim \delta^{-\frac12} \delta^2 |T_S|^{\frac{1}{p}} \sim \delta^{-\frac12} \delta^2 \delta^{-\frac{2}{p}}.
\end{equation}
This shows that the decoupling inequality is false for $p>4$. 

Suppose that the eigenvalues have different signs. Then as in the previous case, we have
\begin{equation}
    \delta  \lesssim \|f\|_{L^p(\mathbb{R}^3)}.
\end{equation}
We need to calculate
\begin{equation}
    \Big( \sum_{S \in \mathcal{S}_{\delta} }\|f_S\|_{L^p}^2 \Big)^{1/2}.
\end{equation}
By the construction of $\mathcal{S}_{\delta}$, we have $|f_S|=\delta^21_{T_S} $
where $T_S$ has dimension $\delta^{A} \times \delta^{1-A} \times \delta$ for some $0 \leq A \leq 1$. Here the number $A$ depends on the choice of $S$. Also, the cardianlity of $\mathcal{S}_{\delta}$ is comparable to $\delta^{-1}$ up to $\delta^{-\epsilon}$ losses. So we have
\begin{equation}
    \Big( \sum_{S \in \mathcal{S}_{\delta} }\|f_S\|_{L^p}^2 \Big)^{1/2} \lesssim_{\epsilon} \delta^{-\epsilon} \delta^{-\frac12} \delta^2 |T_S|^{\frac1p} \sim \delta^{-\epsilon} \delta^{-\frac12}\delta^2 \delta^{-\frac2p}.
\end{equation}
This shows that the decoupling inequality is false for $p>4$, and completes the proof of Item (1). 
\\

Let us move on to Item (2). Since $\phi$ depends only on a variable, after some change of variables and abusing notations, we may assume that $\phi(\xi_1,\xi_2)=\psi(\xi_1)$. Fix $x_2$ and define
\begin{equation}
    g(x_1,x_3):= f(x_1,x_2,x_3).
\end{equation}
Then \eqref{23.10.20} simply follows from \cite[Theorem 1.4]{MR4310157}. Let us show the sharpness. Since $\psi$ is a one-variable polynomial, the zeros of $\psi$ are finite. So we can find an open interval $I \subset [0,1]$ such that for every $\xi_1 \in I$, we have $|\psi''(\xi_1)| \gtrsim 1$. We write $I$ as a union of intervals $J$ with length $\delta^{1/2}$.
Take $\hat{f}$ to be a smooth bump function of the $\delta$-neighborhood of
\begin{equation}
    \{ (\xi_1,\xi_2,\phi(\xi_1,\xi_2) ) :\xi_1 \in I, \xi_2 \in [0,1] \}.
\end{equation} 
Then similar calculations in the proof of Item (1) give the sharpness for $p$. We leave out the details.

Lastly, let us prove Item (3). If $\phi$ does not depend on any variable, then according to our construction, $f_S=f$. So \eqref{23.10.20} is true for $2 \leq p <\infty$. This completes the proof of Proposition \ref{10.23.18}.
\end{proof}

\bibliographystyle{alpha}
\bibliography{reference}

\end{document}